\documentclass[a4paper,10pt]{amsart}
\usepackage{amsmath,amsfonts,amsthm,amssymb}
\usepackage[dvips]{graphicx}
\usepackage{psfrag}

\newtheorem{thm}{Theorem}[section]
\newtheorem{lem}[thm]{Lemma}
\newtheorem{pro}[thm]{Proposition}
\newtheorem{dfn}[thm]{Definition}

\newtheorem{cor}[thm]{Corollary}
\newtheorem{rk}[thm]{Remark}

\def\bb{\begin}
           \def\ea{\end{array}}
          \def\ec{\end{center}}
     \def\ed{\end{description}}
        \def\ee{\end{equation}}
       \def\eea{\end{eqnarray}}
     \def\eeaa{\end{eqnarray*}}
\def\bt{\bb{thebibliography}} \def\et{\end{thebibliography}}
\def\bib{\bibitem}            \def\cal{\mathcal}

\def\Pf {\noindent {\bf Proof}\q}
\def\qed {\hfill $\Box$\vskip5pt}

\def\aaa{\vdash\hskip-2pt\dashv}
\def\h{\pitchfork}
\def\Pf {\noindent {\bf Proof}\quad}
\def\qed {\hfill $\Box$\vskip5pt}

\title{Lyapunov stable chain recurrent classes}
\author{Jiagang Yang}
\date{\today}
\thanks{Partially supported by TWAS-CNPq, FAPERJ}

\address{IMPA, Est. D. Castorina 110, 22460-320 Rio de Janeiro, Brazil}
\email{yangjg\@@impa.br}

\begin{document}

\begin{abstract}
We show that for a $C^1$ residual subset of diffeomorphisms far
away from homoclinic tangency, the stable manifolds of periodic
points cover a dense subset of the ambient manifold. This gives a
partial proof to a conjecture of C. Bonatti.
\end{abstract}

\maketitle


\section{Introduction}\label{s.introduction}
This paper is about generic dynamics, a subject that has been very
active in the last years. The theory of generic dynamics is trying
to give a description of a large class of differential dynamics,
especially it can help us understanding the non-hyperbolic
diffeomorphisms which is one of the most important aim of modern
dynamical theory.

The stable manifold for hyperbolic periodic point is one of the
most basic and important object in differential dynamic, such
submanifold has a special converging property, and the complicated
phenomena: homoclinic intersection just comes from the transverse
intersection between the stable manifold and unstable manifold.
When a diffeomorphism $f$ is hyperbolic, it's well known that the
union of stable manifolds of $f$'s periodic points is dense, but
people discovered that the set of hyperbolic diffeomorphisms are
not dense among differential dynamics, so we want to know that if
the results on the hyperbolic systems can indicate that the same
property will be hold for generic non-hyperbolic systems. Here we
proved that:\\

\noindent{\bf Theorem 1:}{\it There exists a generic subset
$R\subset (\overline{HT^1})^c$ such that for any $f\in R$,
$\bigcup \limits_{p\in Per(f)}W^s(p)$ is dense in
$M$.}\\

These result gives a partial answer to the following Bonatti's
conjecture:\\

\noindent{\bf Conjecture 1} (Bonatti): {\it There exists a generic
subset $R\subset C^1(M)$ such that for any $f\in R$, $\bigcup
\limits_{p\in Per(f)}W^s(p)$ is dense in
$M$.}\\

The Bonatti's conjecture is one step towards the following famous
 conjecture.\\

\noindent{\bf $C^r$ Palis conjecture:} {\it Diffeomorphisms of $M$
exhibiting either a homoclinic tangency or heterdimensional cycle
are $C^r$ dense in the complement of the $C^1$ closure of
hyperbolic systems.}\\

Since until now, almost all the perturbation tools just work in
$C^1$ topology, in this paper we just consider $C^1$
diffeomorphisms and talk about $C^1$ typical phenomena.

In fact, I believe something even stronger than Palis conjecture should be live:\\

\noindent{\bf Conjecture 3 (Tameness conjecture):} {\it There
exists a generic subset $R\subset (\overline{HT^1})^c$ such that
any $f\in R$ is tame.}\\

It's not difficult to get $C^1$ Palis conjecture from tameness
conjecture, but until now we can't prove the tameness conjecture
even in the simplest open set: the small open neighborhood of the
map: {\it linear Anosov map}$|_{T^2}\ \times Id_{S^1}$. In the
flow case, it looks like true in the set ${\mathcal F}^1(M)$.

In the direction of proving the Tameness conjecture, I propose the
following two intermediate problems:\\

\noindent{\bf Conjecture 4:} {\it There exists a generic subset
$R\subset (\overline{HT^1})^c$ such that for any $f\in R$, its
chain recurrent classes are all homoclinic classes.}\\

\noindent{\bf Conjecture 5:} {\it There exists a generic subset
$R\subset (\overline{HT^1})^c$ such that for any $f\in R$, suppose
$C$ is a homoclinic class of $f$, and $i_0=\min \limits_i\{i:\
C\bigcap Per_i(f)\neq \phi\}$, then $C$ has an index $i_0$
dominated splitting $T_CM=E^{s}_{i_0}\oplus E^{cu}_{i_0+1}$ where
$E^{s}_{i_0}$ is contracting.}\\

Here I want to point out that the above two weaker conjectures are
still enough to prove Palis conjecture, now let's show some simple
idea of how to induce $C^1$ Palis conjecture from the above two
conjectures: suppose $f\in R$ and it's far away from
heterdimensional cycle ($f\in (\overline{HC}\bigcup
\overline{HT})^c$), let $C$ be any chain recurrent class of $f$,
then by conjecture 4, $C$ is a homoclinic class, and by $f\in
(\overline{HC})^c$, all the periodic points in $C$ have the same
index $i_0$, then by conjecture 5, $C$ is hyperbolic and has an
index $i_0$ dominated splitting $T_CM=E^{s}_{i_0}\oplus
E^{u}_{i_0+1}$, then it's easy to know $f$ has just finite chain
recurrent classes, so $f$ satisfies Axiom A, $f$ satisfies the
non-cycle condition is just a well known $C^1$ generic result from
\cite{1BC}'s connecting lemma.

The above two conjectures have been proved by \cite{1PS1} when $M$
is a boundless surface (in fact, they proved tameness conjecture
in this case). In higher dimensional manifold they are still far
away to be proved. The following conjectures are weaker more and
look like easier to prove:\\

\noindent{\bf Conjecture 6:} {\it There exists a generic subset
$R\subset (\overline{HT^1})^c$ such that for any $f\in R$, suppose
$C$ is any aperiodic class of $f$, then $C$ has a partial
hyperbolic splitting $T_CM=E^{s}\oplus E^c\oplus E^u$ where
$E^s,E^u\neq \phi$ and $dim(E^c)=1$.}\\

\noindent{\bf Conjecture 7:} {\it There exists a generic subset
$R\subset (\overline{HT^1})^c$ such that for any $f\in R$, suppose
$C$ is a homoclinic class of $f$ and $i_0=\min \limits_i\{i:\
C\bigcap Per_i(f)\neq \phi\}$, then $C$ has an index $i_0$
dominated splitting $T_CM=E^{cs}_{i_0}\oplus E^{cu}_{i_0+1}$, and
either $E^{cs}_{i_0}$ is contracting or $E^{cs}_{i_0}$ has a
codimension-1 sub-dominated splitting
$E^{cs}_{i_0}=E^{s}_{i_0-1}\oplus E^{c}_1$ where $E^{s}_{i_0-1}$
is hyperbolic
and $dim(E^c_1|_C)=1$}.\\

All the conjectures above just talk about general chain recurrent
classes, before we prove them, we should check them in some
special situation. In this paper we'll use a special chain
recurrent class: Lyapunov stable chain recurrent class to check
these conjectures, and we can show that for this special kind of
chain recurrent class Conjecture 4 and half of Conjecture 7 are
right, they give some evidence that the above conjectures may be
right.
The precisely statements are following:\\

\noindent{\bf Theorem 2:} {\it There exists a generic subset
$R\subset (\overline{HT^1})^c$ such that for any $f\in R$, its
Lyapunov stable chain recurrent classes should be homoclinic
classes.}\\

\noindent{\bf Theorem 3:} {\it There exists a generic subset
$R\subset (\overline{HT^1})^c$ such that for any $f\in R$, suppose
$C$ is any Lyapunov stable homoclinic class of $f$, let $i_0=\min
\limits_i\{i:\ C\bigcap Per_i(f)\neq \phi\}$, then $C$ has an
index $i_0$ dominated splitting $T_CM=E^{cs}_{i_0}\oplus
E^{cu}_{i_0+1}$, and
\begin{itemize}
\item either $E^{cs}_{i_0}$ is contracting and $C$ is an index
$i_0$ fundamental limit

\item or $E^{cs}_{i_0}$ has a codimension-1 sub-dominated
splitting $E^{cs}_{i_0}=E^{s}_{i_0-1}\oplus E^{c}_1$ where
$E^{s}_{i_0-1}$ is contracting and $dim(E^c_1|_C)=1$, $C$ is an
index $i_0-1$ and index $i_0$ fundamental limit.
\end{itemize}}

In $\S$3 we'll state some generic properties and give an important
technique lemma, its proof will be given in $\S$7. In $\S$4 I'll
introduce some properties for fundamental limit and Crovisier's
central model, in $\S$5 I'll state the main lemma and use it to
prove theorem 1,2,3. The
proof of the main lemma is given in $\S$6. \\

\noindent{\bf Acknowledgements:} This paper is part of the
author's thesis. I would like to thank my advisor Professor
Marcelo Viana for his support and enormous encouragements during
the preparation of this work. I would like to thank Professor
Shaobo Gan for checking the details of the proof and finding out
an essential gap in the original argument which is crucial to the
work. I also thank Professors Jacob Palis, Lan Wen, Enrique R.
Pujals, Lorenzo Diaz, Christian Bonatti for very helpful remarks.
Finally I wish to thank my wife, Wenyan Zhong, for her help and
encouragement.

\section{Definitions and Notations}

Let $M$ be a compact boundless Riemannian manifold, since when $M$
is a surface \cite{1PS1} has proved that hyperbolic
diffeomorphisms are open and dense in $C^1(M)\setminus
\overline{HT}$, we suppose $dim(M)=d>2$ in this paper. Let
$Per(f)$ denote the set of periodic points of $f$ and $\Omega (f)$
the non-wondering set of $f$, for $p \in Per(f)$, $\pi (p)$ means
the period of $p$. If $p$ is a hyperbolic periodic point, the
index of $p$ is the dimension of the stable bundle. We denote
$Per_i(f)$ the set of the index $i$ periodic points of $f$, and we
call a point $x$ is an index $i$ preperiodic point of $f$ if there
exists a family of diffeomorphisms
$g_n\stackrel{C^1}{\longrightarrow}f$, where $g_n$ has an index
$i$ periodic point $p_n$ and ${p_n}\longrightarrow x$. $P^*_i(f)$
is the set of index $i$ preperiodic points of $f$.

\begin{rk}\label{2.1}It's easy to know $\overline{P_i(f)}\subset P^*_i(f)$.
\end{rk}

Let $\Lambda$ be an invariant compact set of $f$, we call
$\Lambda$ is an index $i$ fundamental limit if there exists a
family of diffeomorphisms ${g_n}$ $C^1$ converging to $f$, $p_n$
is an index $i$ periodic point of $g_n$ and $Orb(p_n)$ converge to
$\Lambda$ in \emph{Hausdorff} topology. So if $\Lambda (f)$ is an
index $i$ fundamental limit, we have $\Lambda (f)\subset
P^*_i(f)$. $\Lambda$ is a minimal index $i$ fundamental limit if
$\Lambda (f)$ is an index $i$ fundamental limit and any invariant
compact subset $\Lambda_0\varsubsetneq\Lambda$ is not an index $i$
fundamental limit. In \cite{1Y2} we have showed the following
result:

\begin{lem}\label{2.2}
Any index $i$ fundamental limit contains a minimal index $i$
fundamental limit.
\end{lem}

For two points $x,y\in M$ and some $\delta >0$, we say there
exists a $\delta$-pseudo orbit connects $x$ and $y$ if there exist
points $x=x_0,x_1,\cdots,x_n=y$ such that
$d(f(x_i),x_{i+1})<\delta$ for $i=0,1,\cdots ,n-1$, and we denote
it $x \underset{\delta}{\dashv}y$. We say $x\dashv y$ if for any
$\delta >0$ we have $x \underset{\delta}{\dashv} y$ and denote $x
\aaa y$ if $x\dashv y$ and $y\dashv x$. A point $x$ is called a
chain recurrent point if $x\aaa x$. $CR(f)$ denotes the set of
chain recurrent points of $f$, it's easy to know that $\aaa$ is a
closed equivalent relation on $CR(f)$, and every equivalent class
of such relation should be compact and called chain recurrent
class. A chain recurrent class $C$ of $f$ is called Lyapunov
stable if there exists a family of neighborhoods $\{U_n\}$ of $C$
satisfying: \begin{itemize} \item[a)] $\overline{U_{n+1}}\subset
U_n$,

\item[b)]$\bigcap U_n=C$,

\item[c)]$f(\overline{U_n})\subset U_n$.
\end{itemize}

\begin{rk}\label{2.3}
Conley proved that any homeomorphism $f$ has at least one Lyapunov
stable chain recurrent class.
\end{rk}

\begin{lem} \label{Lya}

Let $C$ be a Lyapunov stable chain recurrent class of $f$, then if
$y\in W^u(C)$ (that means $\lim \limits_{i\rightarrow
\infty}\min\limits_{z\in C}\{d(f^{-i}(y),z)\}\longrightarrow 0$),
we have $y\in C$.
\end{lem}
\Pf: For any $U_n$ the neighborhood of $C$ given in the definition
of Lyapunov stable chain recurrent class, there exists an $i>0$
such that $f^{-i}(y)\in U_n$, then $y\in f^i(U_n)\subset U_n$, so
$y\in \bigcap\limits_n U_n=C$. \qed

Let $K$ be a compact invariant set of $f$, and $x,y$ are two
points in $K$, we denote $x \underset{K}{\dashv} y$ if for any
$\delta
>0$, we have a $\delta$ -pseudo orbit in $K$ connects $x$ and $y$.
If for any two points $x,y\in K$ we have $x \underset{K}{\dashv}
y$, we call $K$ a chain recurrent set. Let $C$ be a chain
recurrent class of $f$, we say $C$ is an aperiodic class if $C$
does not contain periodic point.

Let $\Lambda$ be an invariant compact set of $f$, for $l\in
\mathbb{N}$, $0<\lambda <1$ and $1\leq i<d$, we say $\Lambda$ has
an index $i-(l,\lambda)$ dominated splitting if we have a
continuous invariant splitting $T_{\Lambda}M=E\oplus F$ where
$dim(E_x)=i$ for any $x\in \Lambda$ and $\parallel
Df^l|_{E(x)}\parallel\cdot
\parallel Df^{-l}|_{F(f^{l}x)}\parallel <\lambda$ for all $x\in \Lambda$.
For simplicity, sometimes we just say $\Lambda (f)$ has an index
$i$ dominated splitting. A compact invariant set can have many
dominated splittings, but for fixed $i$, the index $i$ dominated
splitting is unique.

  We say a diffeomorphism $f$ has $C^r$ tangency if $f\in C^r(M)$, $f$ has hyperbolic
  periodic point $p$ and there exists a non-transverse intersection between
  $W^s(p)$ and $W^u(p)$. $HT^r$ denote the set of the diffeomorphisms which have
  $C^r$ tangency, usually we just use $HT$ denote $HT^1$. We call a
diffeomorphism $f$
  is far away from tangency if $f\in C^1(M)\setminus \overline{HT}$. The following proposition shows the
  relation between dominated splitting and far away from
  tangency.

\begin{pro}\label{2.4} (\cite{1W1}) $f$ is $C^1$ far away from tangency if
and only if there exists $(l,\lambda)$ such that $P^*_i(f)$ has
index $i-(l,\lambda)$ dominated splitting for $0<i<d$.\end{pro}

  Usually dominated splitting is not a hyperbolic splitting, Ma\~{n}\'{e} showed
  that in some special case, one bundle of the dominated splitting is
  hyperbolic.

\begin{pro}\label{2.5}(\cite{1M3}) Suppose $\Lambda (f)$ has an index $i$
dominated splitting $E\oplus F\  (i\neq 0)$, if $\Lambda
(f)\bigcap P^*_j(f)=\phi$ for $0\leq j<i$, then $E$ is a
contracting bundle.\end{pro}

\section{Generic properties}
Here we'll introduce some $C^1$ generic properties.

 For a topology space $X$, we call a set $R\subset X$ is a generic subset
  of $X$ if $R$ is countable intersection of open and dense subsets of $X$,
  and we call a property is a generic property of $X$ if there exists
  some generic subset $R$ of $X$ holds such property. Especially, when
  $X=C^1(M)$ and $R$ is a generic subset of $C^1(M)$, we just call $R$ is
  $C^1$
  generic, and we call any generic property of $C^1(M)$ 'a $C^1$ generic
  property' or 'the property is $C^1$ generic'.

  It's easy to know that if $R$ is $C^1$ generic and $R_1$ is a generic subset of $R$, then $R_1$ is also $C^1$ generic.

  At first let's state some well known $C^1$ generic properties.

\begin{pro}\label{3.1} There is a $C^1$ generic subset $R_0$ such
that for any $f\in R_0$, one has
\begin{itemize}
  \item[1)]  $f$ is Kupka-Smale (every periodic
point $p$ in $Per(f)$ is hyperbolic and the invariant manifolds of
periodic points are everywhere transverse).
  \item[2)] $CR(f)=\Omega =\overline{Per(f)}$.
  \item[3)] $P^*_i(f)=\overline{P_i(f)}$
  \item[4)]  any chain recurrent set is the Hausdorff
limit of periodic orbits.
  \item[5)]  any index $i$ fundamental limit is the
Hausdorff limit of index $i$ periodic orbits of $f$.
  \item[6)] any chain
recurrent class containing a periodic point $p$ is the homoclinic
class $H(p,f)$.
  \item[7)] suppose $C$ is a homoclinic class of $f$, and $i_0=min\{i:C\bigcap Per_i(f)\neq
  \phi\},\ i_1=max\{i: C\bigcap Per_i(f)\neq \phi\},$
then for any $i_0\leq i\leq i_1$, we have $C\bigcap Per_i(f)\neq
\phi$ and $C$ is index $i$ fundamental limit.
  \item[8)] if all the Lyapunov stable chain recurrent classes of $f$
are homoclinic classes, then $\bigcup \limits_{p\in Per(f)}W^s(p)$
is dense in $M$.
\end{itemize}
\end{pro}
\Pf 1) comes from Kupka-Smale theorem, 2) is proved in \cite{1BC},
3),4),5),6) are all well known, 7) is proved in \cite{1ABCD}, 8)
is proved in \cite{1MP}. \qed
  By proposition \ref{3.1}, for any $f$ in $R_0$, every chain recurrent class
  $C$
  of $f$ is either an aperiodic class or a homoclinic class. If $\# (C)=\infty$,
  we say $C$ is non-trivial.

  The following technique lemma gives a new $C^1$ generic property
  whose proof would be given in $\S$7.
  \begin{lem}\label{3.2} (Technique lemma). There exists a generic subset $R^\prime_0$ of
$R_0$ such that for $f\in R^\prime_0$, suppose $C$ is a
non-trivial chain recurrent class of $f$, $\Lambda \varsubsetneq
C$ is a compact chain recurrent set without periodic point, then
for $0<s<1$ and any point $y\in (C\setminus \Lambda)\bigcap
W^{s(u)}(\Lambda)$, for any small neighborhood $O$ of $y$ and any
small neighborhood $V$ of $\Lambda$, there exists a periodic point
$q$ of $f$ satisfying $Orb(q)\bigcap O\neq \phi$, and $\frac{\#{\{
Orb(q)\bigcap
   V\}}}{\pi(q)}>s$.\end{lem}

Since $R^\prime_0$ is a generic subset of $R_0$ and $R_0$ is $C^1$
generic, $R^\prime_0$ is a $C^1$ generic subset also.

\begin{cor}\label{3.3}
There exist a generic subset $R\subset R_0^\prime \setminus
\overline{HT}$ such that for $f\in R$, if $C$ is a chain recurrent
class of $f$, $\Lambda\varsubsetneq C$ is a non-trivial minimal
set with partial hyperbolic splitting $E^s_i\oplus E^c_1\oplus
E^u_{i+1}$ where $dim(E^c_1(\Lambda))=1$ and $E^c_1(\Lambda)$ is
not hyperbolic, then $W^u(\Lambda)\bigcap C\subset P^*_i\bigcap
P^*_{i+1}$ and
\begin{itemize}
\item either $C$ contains index $i+1$ or index $i$ periodic point
and it's an index $i$ fundamental limit,

\item or for any $y\in (C\bigcap W^u(\Lambda))\setminus \Lambda$,
and $\{V_n\}$ is a family of neighborhoods of $\Lambda$ satisfying
$\overline{V_{n+1}}\subset V_n$ and $\bigcap V_n=\Lambda$, there
exists $\{q_n\}$ a family of index $i$ (or $i+1$) periodic points
of $f$ such that $(y\bigcup \Lambda)\subset \lim
\limits_{n\rightarrow \infty}Orb(q_n)$ and $\lim
\limits_{n\rightarrow \infty}\frac{\#\{Orb(q_n)\bigcap
V_n\}}{\pi(q_n)}\longrightarrow 1^-$.
\end{itemize}
\end{cor}
\Pf: At first let's suppose $f\in R_0^\prime \setminus
\overline{HT}$. When $i=0$ (or $i+1=d$), theorem 1 of \cite{1Y2}
has shown $C$ contains index 1 ($d-1$) periodic point and $C$ is
an index 0 and index 1 (index $d$ and index $d-1$) fundamental
limit, so from now we suppose $E^s_i|_\Lambda,
E^u_{i+2}|_\Lambda\neq \phi$, and here we just prove the above
result for case $i$, the proof of the case $i+1$ is similar.

Fix any $y\in C\bigcap W^u(\Lambda)\setminus \Lambda$ and $V_n$ is
a family of neighborhood of $\Lambda$ such that
$\overline{V_{n+1}}\subset V_n$ and $\bigcap \limits_{n\geq
1}V_n=\Lambda$, choose $\varepsilon_n>0$ and $0<s_n<1$ satisfying
$\varepsilon_n\longrightarrow0^+$ and $s_n\longrightarrow1^-$. By
the technique lemma, there exists a family of periodic points
$\{q_n(f)\}$ such that $y\bigcup \Lambda \subset \lim
\limits_{n\rightarrow \infty}Orb(q_n)$ and $\{q_n\}$ satisfies
$\frac{\#\{Orb(q_n)\bigcap V_n\}}{\pi(q_n)}>s_n$. We can let all
the $q_n(f)$ have the same index $j$, we suppose $j\geq i$, since
the proof of the other case is the same.

 Let $j_1=\min \limits_{j\geq i}\{j:$ there exists a family of
$C^1$ diffeomorphism $g_n$ such that $\lim \limits_{n\rightarrow
\infty}g_n\longrightarrow f$ and $g_n$ has an index $j$ periodic
point $p_n(g_n)$ such that $\lim \limits_{n\rightarrow
\infty}Orb_{g_n}(q_n(g_n))\supset y\bigcup \Lambda$ and
$\frac{\#\{Orb_{g_n}(q_n)\bigcap V_n\}}{\pi(q_n)}>s_n\}$.

We claim that
\begin{itemize}
\item[(a)] either $C$ contains index $i+1$ or index $i$ periodic
point and it's an index $i$ fundamental limit,

\item[(b)] or $j_1=i$.
\end{itemize}
\noindent{\bf Proof of the claim}

\noindent{$\bullet$}  If $j_1=i$, we get (b).\\
\noindent{$\bullet$} If $j_1>i$, we'll show (a) is true.

Suppose $g_n$ is the family of diffeomorphisms and $q_n(g_n)$ is
the index $j_1$ periodic point of $g_n$ given in the definition of
$j_1$. Let $\lim \limits_{n\rightarrow \infty}Orb_{g_n}(q_n)=C_0$,
then $C_0\subset P_{j_1}^*$, by proposition \ref{2.4}, $C_0$ has
an index $j_1$ dominated splitting $E^{cs}_{j_1}\oplus
E^{cu}_{j_1+1}|_{C_0}$.

 By the definition of $j_1$ and
Franks lemma, we know that
$\{Dg_n|_{E^{cs}_{j_1}(Orb_{g_n}(q_n))}\}^\infty_{n=1}$ is stable
contracting. By lemma 4.9, lemma 4.10 and remark 4.11 of
\cite{1Y2}, there exist $N_0$, $l$, $0<\lambda<1$ such that for
$\pi_{g_n}(q_n)>N_0$, we have $q_n^\prime\in Orb_{g_n}(q_n)$
satisfying $\prod \limits^{s-1}_{j=0}\Vert
Dg^l_n|_{E^{cs}_{j_1}(g^{jl}_n(q_n^\prime))}\Vert \leq \lambda^s$
for $s\geq 1$. Since $\Lambda$ is minimal and non-trivial, from
$\Lambda\subset \lim \limits_{n\rightarrow
\infty}Orb_{g_n}(q_n(g_n))$, we know $\lim
\limits_{n\rightarrow\infty}\pi_{g_n}(q_n(g_n))\longrightarrow
\infty$, so we can suppose $\pi_{g_n}(q_n(g_n))>N_0$ always. The
above point $q_n^\prime$ is called hyperbolic time for bundle
$E^{cs}_{j_1}$, its existence comes from Pliss lemma, since
$Orb_{g_n}(q_n)$ stays a lot of time in $V_n$, so in fact from the
Pliss lemma we can always choose $q_n^\prime\in V_n$, then we can
suppose $\lim \limits_{n\rightarrow \infty}q_n^\prime=x_0\in
\Lambda$, by $\lim \limits_{n\rightarrow
\infty}g_n\overset{C^1}{\longrightarrow}f$, we have
\begin{equation} \label{17}
\prod \limits^{s-1}_{j=0}\Vert
Df^l|_{E^{cs}_{j_1}(f^{jl}(x_0))}\Vert \leq \lambda^s\ \text{for}\
s\geq 1. \end{equation}

Since $\Lambda$ has two dominated splitting $(E^{cs}_{i}\oplus
E^c_1)\oplus E^{cu}_{i+2}$ and $E^{cs}_{j_1}\oplus E^{cu}_{j_1+1}$
with $j_1\geq i+1$, by lemma 4.30 of \cite{1Y2}, we have that
$E^{cs}_i\oplus E^c_1 \subset E^{cs}_{j_1}$, so by \eqref{17}, we
get $\prod \limits^{s-1}_{j=0}\Vert Df^l|_{(E^{cs}_{i_0}\oplus
E^c_1)(f^{jl}(x_0))}\Vert \leq \lambda^s$ for $s\geq 1$. By
$\Lambda$ is minimal and $E^c_1|_\Lambda$ is not hyperbolic, the
splitting $(E^{cs}_{i}\oplus E^c_1)\oplus E^{cu}_{i+2}|_\Lambda$
satisfies all the assumptions of weakly selecting lemma, by weakly
selecting lemma given in \cite{1Y2} and corollory 4.26 there, $C$
contains index $i+1$ periodic point and $C$ is an index $i$
fundamental limit, so $C$ satisfies (a). \qed

Now with a generic argument like we'll do in $\S$7.1, in the proof
of above claim we can replace $R_0^\prime\setminus \overline{HT}$
by a generic subset $R\subset R_0^\prime\setminus \overline{HT}$
such that if $f\in R$ and (a) is false, $f$ itself will have a
family of index $i$ periodic points $\{q_n\}$ such that $(y\bigcup
\Lambda)\subset \lim \limits_{n\rightarrow \infty}Orb(q_n)$ and
$\lim \limits_{n\rightarrow \infty}\frac{\#\{Orb(q_n)\bigcap
V_n\}}{\pi(q_n)}\longrightarrow 1^-$.\qed

We'll show the generic set $R$ satisfies theorem 1, 2 and 3.

\section{Fundamental limit and Crovisier's central model}
\subsection{The minimal index $j_0$ fundamental limit}

Let $f\in R$, $C$ is any non-trivial chain recurrent class of $f$,
suppose $j_0=\min \limits_j \{j:\ C\bigcap P^*_j\neq \phi\}$ and
$\Lambda$ be a minimal index $j_0$ fundamental limit, by lemma
\ref{2.2}, such set always exists. Now we'll recall some results
about $j_0$ and the set $\Lambda$, they are all given in
\cite{1Y2}

\begin{lem}\label{4.1} Suppose $f\in R$, $C$ is a chain recurrent
class of $f$, $j_0=\min \limits_j \{j:\ C\bigcap P^*_j\neq
\phi\}$, $\Lambda$ is a minimal index $j_0$ fundamental limit in
$C$, then
\begin{itemize}
\item either $\Lambda$ is a non-trivial minimal set with partial
hyperbolic splitting $E^s_{j_0}\oplus E^c_1\oplus E^u_{j_0+2}$

\item or $C$ contains a periodic point with index $j_0$ or $j_0+1$
and $C$ is an index $j_0$ fundamental limit.
\end{itemize}

\end{lem}

\begin{lem}\label{4.3}
Suppose $f\in R$, $C$ is a non-trivial chain recurrent class of
$f$, if $C\bigcap P^*_0\neq \phi$, then $C$ should be a homoclinic
class containing index 1 periodic points and $C$ is an index 0
fundamental limit.
\end{lem}

\subsection{Partial hyperbolic splitting and Crovisier's central model}

$f\in C^1(M)$, Suppose $\Lambda$ is a minimal set of $f$ with
partial hyperbolic splitting $E^s_{j_0}\oplus E^c_1\oplus
E^u_{j_0+2}$ where $E^s_{j_0},E^u_{j_0+2}\neq \phi$,
$dim(E^c_1|_\Lambda)=1$ and $E^c_1(\Lambda)$ is not hyperbolic,
let $C$ be the chain recurrent class containing $\Lambda$ and
$V_0$ be a small neighborhood of $\Lambda$, then the maximal
invariant set of $\overline{V_0}$: $\Lambda_0=\bigcap \limits_j
f^j(\overline{V_0})$ will have a partial hyperbolic splitting
$E^s_{j_0}\oplus E^c_1\oplus E^u_{j_0+2}$ also. In fact, we can
extend such splitting to $\overline{V_0}$ (it's not invariant
anymore). For every point $x\in \overline{V_0}$, we define some
cones on its tangent space $C^i_a(x)=\{v|v\in T_xM,$ there exists
$v^\prime \in E^i(x)$ such that
$d(\frac{v}{|v|},\frac{v^\prime}{|v^\prime|})<a\}_{i=s,c,u,cs,cu}$.
When $a$ is small enough, $C^i_a(x)\bigcap C^j_a(x)=\phi\ _{(i\neq
j= s, c, u)}$, $C^{cs}_a(x)\bigcap C^u_a(x)=\phi$,
$C^{cu}_a(x)\bigcap C^s_a(x)=\phi$ for any $x\in \overline{V_0}$,
and $Df(C^i_a(x))\subset C^i_a(f(x))\ _{i=u,cu}$,
$Df^{-1}(C^i_a(x))\subset C^i_a(f^{-1}(x))\ _{i=s,cs}$ for $x\in
\Lambda_0$.

We say a submanifold $D^i$ ($i=s,c,u,cs,cu$) tangents with cone
$C^i_a$ when $dim D^i=dim(E^i)$ (we denote $E^{cs}=E^c_1\oplus
E^s$, $E^{cu}=E^c_1\oplus E^u$) and for $x\in D^i$, $T_xD^i\subset
C^i_a(x)$. For simplicity, sometimes we just call it $i$-disk,
especially when $i=c$, we call $D^c$ a central curve. We say an
$i$-disk $D^i$ has center $x$ with size $\delta$ if $x\in D^i$,
and respecting the Riemannian metric restricting on $D^i$, the
ball centered on $x$ with radius $\delta$ is in $D^i$. We say an
$i$-disk $D^i$ has center $x$ with radius $\delta$ if $x\in D^i$,
and respecting the Riemannian metric restricting on $D^i$, the
distance between any point $y\in D^i$ and $x$ is smaller than
$\delta$.

We say a smooth central curve $\gamma$ is a central segment if
$f^i(\gamma)\subset V_0$ and $f^i(\gamma)$ is a central curve for
any $i\in  \mathbb{Z}$, so if $\gamma$ is a central segment,
$\gamma\subset \Lambda_0$, and it's easy to know
$T_x\gamma=E^c_1(x)$ for any $x\in \gamma$. We say a smooth
central curve $\gamma$ is a positive(negative) central segment if
$f^i(\gamma)\subset V_0$ and $f^i(\gamma)$ is a central curve for
any $i\geq(\leq) 0$, so if $\gamma$ is a positive (negative)
central segment, $\gamma\subset \bigcap \limits_{-\infty}^0
f^i(\overline{V_1})\ (\bigcap \limits_0^\infty
f^i(\overline{V_1}))$,

Now let's consider the orientation of the central bundle
$E^c_1(\Lambda)$.

\begin{dfn}\label{4.4}
We say $E^c_1(\Lambda)$ has an $f$-orientation if $E^c_1(\Lambda)$
is orientable and $Df$ preserves its orientation.

\end{dfn}
\begin{lem}\label{4.5}
For a compact neighborhood $V_1$ of $\Lambda$ satisfying
$V_1\subset V_0$ and let $\Lambda_1=\bigcap
\limits_{i=-\infty}^\infty f^i(V_1)$, $\Lambda_1^+=\bigcap
\limits_{i=-\infty}^0 f^i(V_1)$, $\Lambda_1^-=\bigcap
\limits_{i=0}^\infty f^i(V_1)$, then there exist $\delta_0>0$,
$\delta_0/2>\delta_1>\delta_2>0$ such that they satisfy the
following properties:
\begin{itemize}
\item[a)] If $E^c_1(\Lambda)$ has an $f$ orientation,
$E^c_1(\Lambda_1)$ has an $f$ orientation also.

\item[b)]for any $x\in V_1$, $B_{\delta_0}(x)\subset V_0$ and
$E^c_1(B_{\delta_0}(x))$ is orientable, so it gives orientation
for any central curve in $B_{\delta_0}(x)$, and we suppose
$\delta_0$ is small enough such that any central curve in
$B_{\delta_0}(x)$ never intersects with itself.

\item[c)] for any $x\in \Lambda_1^+$, $x$ has $\delta_1$ uniform
size of strong stable manifold $W^{ss}_{\delta_1}(x)$ and
$W^{ss}_{\delta_1}(x)$ is an $s$ disk; for any $x\in \Lambda_1^-$,
$x$ has $\delta_1$ uniform size of strong unstable manifold
$W^{uu}_{\delta_1}(x)$ and $W^{uu}_{\delta_1}(x)$ is an $u$ disk.

\item[d)] for any $x\in \Lambda_1$, there exists a central curve
$l_{\delta_1}(x)$ with center $x$ and radius $\delta_1$, such that
there exists a continuous function
$\Phi^c:\Lambda_1\longrightarrow Emb^1(I,M)$ satisfying
$\Phi^c(x)=l_{\delta_1}(x)$ where $x\in \Lambda_1$, and if let
$l_{\delta_2}(x)\subset l_{\delta_1}(x)$ be the central curve with
center $x$ and radius $\delta_1$, then $f(l_{\delta_2}(x))\subset
l_{\delta_1}(f(x))$ and $f^{-1}(l_{\delta_2}(x))\subset
l_{\delta_1}(f^{-1}(x))$.

\item[e)]For any $0<\varepsilon<\delta_1$, there exists $\delta>0$
such that for any positive central segment $\gamma\subset
\Lambda_1^+$ with $\varepsilon<length(\gamma)<\delta_1$,
$W^s_{loc}(\gamma)=\bigcup \limits_{x\in
\gamma}W^{ss}_{\delta_1}(x)$ is a $cs$ disk with uniform size
$\delta$, and for any $x\in Int(\gamma)$, there exists
$\delta_x>0$ such that for any $y\in B_{\delta_x}(x)\bigcap
\Lambda_1$, we have $W^{uu}_{\delta_1}(y)\h W^s_{loc}(\gamma)\neq
\phi$. And if $C_0$ is a invariant compact subset containing
$\Lambda$ and has the following dominated splitting
$E^s_{j_0}\oplus E^c_1\oplus W^u_{j_0+2}$, then there exists $U_n$
a small neighborhood of $C_0$ such that any $z\in C^\prime=\bigcap
\limits_{i}f^i(U_n)$ will have uniform size of strong unstable
manifold $W^{uu}_{\delta_1}(z)$ and if $z\in
B_{\delta_x}(x)\bigcap C^\prime$, we still have
$W^{uu}_{\delta_1}(z)\h W^s_{loc}(\gamma)\neq \phi$ (If
$\gamma\subset \Lambda^-$, we'll have $W^{ss}_{\delta_1}(z)\h
W^u_{loc}(\gamma)\neq \phi$).
\end{itemize}
\end{lem}

\Pf
 a), b) are obviously, c) is \cite{1HPS}'s result about strong stable
 manifold theorem, d) is \cite{1HPS}'s result about cental
 manifolds, the first part of e) is the stable manifold theorem for normally hyperbolic
 submanifold; about the second part, when $U_n$ is small enough, $C_0^\prime$ will have the dominated
 splitting $E^s_{j_0}\oplus E^c_1\oplus
W^u_{j_0+2}$ also, and we can even extend such splitting to $U_n$
and get two cones $C^u_{a_0}|_{C_0^\prime}$ and
$C^{cs}_{a_0}|_{C_0^\prime}$ which match the respectively cones in
$V_1$, then when $\delta_1$ is small enough, any $z\in C_0^\prime$
will have uniform size of strong unstable manifold
$W^{uu}_{\delta_1}(z)$ and it's an u-disk (tangents the cone
$C^u_{a_0}|_{U_n}$), so when $z\in C_0^\prime$ near $x$ enough,
we'll have $W^{uu}_{\delta_1}(z)\h W^s_{loc}(\gamma)\neq \phi$.
\qed

 Now let's introduce Crovisier's result, we divide the statement
 to two cases: $E^c_1(\Lambda)$ has an $f$ orientation or not. At
 first, suppose $E^c_1(\Lambda)$ has an $f$ orientation, and we
 call the direction right.

\begin{lem}\label{4.6}($E^c_1(\Lambda)$ has an $f$ orientation):
$0<\delta_2<\delta_1<\delta_0/2$ are given by lemma \ref{4.5}, and
$l_{\delta_1}(x)\ _{(x\in \Lambda_1)}$ are given there also, let
$l^+_{\delta_1}(x)\subset l_{\delta_1}(x)$ be the central curve in
the right of $x$, then
\begin{itemize}
\item[a)] either for some $x_0\in \Lambda$, there exists a central
segment $\gamma_0\subset l^+_{\delta_1}(x_0)$ where $\gamma_0$
contains $x_0$ and $\gamma_0\subset \Lambda_0$, in fact,
$\gamma_0$ is in the same chain recurrent class with $\Lambda$
respect the map $f|_{V_0}$.

\item[b)] or for any $x\in \Lambda_1$ there exists a central curve
$\gamma_x^{+}\subset l^+_{\delta_1}(x)$ such that $\gamma_x^{+}$
containing $x$, $\gamma_x^{+}\ _{(x\in \Lambda_1)}$ is a family of
smooth curve and they are $C^0$ continuously depend on $x\in
\Lambda_1$, and either $f(\overline{\gamma_x^{+}})\subset
\gamma_{f(x)}^{+}$ for all $x\in \Lambda_1$ or
$f^{-1}(\overline{\gamma_x^{+}})\subset \gamma_{f^{-1}(x)}^{+}$
for all $x\in \Lambda_1$.
\end{itemize}
\end{lem}

In the case b) of lemma \ref{4.6}, if we have
$f(\overline{\gamma_x^{+}})\subset \gamma_{f(x)}^{+}$, we call the
right central curve is 1-step contracting, if
$f^{-1}(\overline{\gamma_x^{+}})\subset \gamma_{f^{-1}(x)}^{+}$,
we call it's one step expanding.

\begin{lem}\label{4.7}(\cite{1C2},\cite{1Y2})
When $f\in R$, and a) of lemma \ref{4.6} happens, then $C$ is a
homoclinic class containing index $i_0$ or $i_0+1$ periodic point
and $C$ is an index $i_0$ fundamental limit.
\end{lem}

\begin{lem}\label{4.8}($E^c_1(\Lambda)$ has non $f$-orientation)
$0<\delta_2<\delta_1<\delta_0/2$ are given by lemma \ref{4.5}, and
$l_{\delta_1}(x)\ _{(x\in \Lambda_1)}$ are given there also, then
\begin{itemize}
\item[a)]either for some $x_0\in \Lambda$, there exists a central
segment $\gamma_0\subset l_{\delta_1}(x_0)$ such that $x_0\in
\gamma_0$ and $\gamma_0\subset \Lambda_1$, and if $f\in R$, then
$C$ is a homoclinic class containing index $i_0$ or $i_0+1$
periodic point and $C$ is an index $i_0$ fundamental limit.

\item[b)] or for every $x\in \Lambda_1$ there exists a central
curve $\gamma_x^{}\subset l_{\delta_1}(x)$ containing $x$ and
$\gamma_x\ _{(x\in \Lambda_1)}$ is a family of smooth curve $C^0$
continuously depend on $x\in \Lambda_1$, and either
$f(\overline{\gamma_x})\subset \gamma_{f(x)}$ for all $x\in
\Lambda_1$ or $f^{-1}(\overline{\gamma_x})\subset
\gamma_{f^{-1}(x)}$ for all $x\in \Lambda_1$.
\end{itemize}
\end{lem}

\section{Proof of theorem 1, 2 and 3}
At first, let's state the main lemma, its proof would be given in
$\S$6.

\begin{lem}\label{5.1}(The main lemma)
Suppose $f\in R$, $C$ is a non-trivial Lyapunov stable chain
recurrent class of $f$, let $j_0=\min \limits_j\{j:C\bigcap
P^*_j\neq \phi\}$, then $C$ contains index $j_0$ or $j_0+1$
periodic point and $C$ is an index $j_0$ fundamental limit.
\end{lem}

It's easy to see that theorem 2 is a simply corollary of the main
lemma.

Now I'll show the proof of theorem 1:

\noindent{\bf Proof of Theorem 1:} It's just a corollary of
generic property 8) of proposition \ref{3.1} and theorem 2. \qed

 \noindent{\bf Proof of Theorem 3:} Recall $j_0=\min
\limits_j\{j:C\bigcap P^*_j\neq \phi\}$ and $i_0=\min
\limits_{i}\{i:C\bigcap Per_i(f)\neq \phi\}$, so $j_0\leq i_0$, by
lemma \ref{5.1}, we have $j_0\geq i_0-1$.

So either $j_0=i_0$ or $j_0=i_0-1$.

When $j_0=i_0$, then by generic property 6) of proposition
\ref{3.1}, $C\subset \overline{Per_{i_0}(f)}\subset P^*_{i_0}(f)$.
By proposition \ref{2.4} and $f\in R\subset (\overline{HT})^c$,
$C$ has an index $i_0$ partial hyperbolic splitting
$T_CM=E^{cs}_{i_0}\oplus E^{cu}_{i_0+1}$. By the definition of
$j_0$ and the assumption $i_0=j_0$, we know $C\bigcap P^*_j=\phi$
for $j<i_0$, so from proposition \ref{2.5}, $E^{cs}_{j_0}|_C$ is
hyperbolic, we denote it by $E^s_{i_0}|_C$, then on $C$ we have
the following dominated splitting $T_CM=E^{s}_{i_0}\oplus
E^{cu}_{i_0+1}$. And since $C$ contains index $i_0$ periodic
point, $C$ is an index $i_0$ fundamental limit.

When $j_0=i_0-1$, by lemma \ref{5.1}, $C$ is an index $i_0-1$
fundamental limit, so $C\subset P^*_{i_0-1}$ and we've known that
$C$ contains index $i_0$ periodic point, so $C\subset P^*_{i_0}$,
then $C\subset P^*_{i_o-1}\bigcap P^*_{i_0}$, from $f\in R\subset
(\overline{HT})^c$ and proposition \ref{2.4}, $C$ has an index
$i_0-1$ dominated splitting $E^{cs}_{i_o-1}\oplus E^{cu}_{i_0}|_C$
and an index $i_0$ dominated splitting $E^{cs}_{i_0}\oplus
E^{cu}_{i_0+1}|_C$. Let $E^c_1|_C=E^{cs}_{i_0}\bigcap
E^{cu}_{i_0}|_C$, then $C$ will have the following dominated
splitting $E^{s}_{i_0-1}\oplus E^{c}_{1}\oplus E^{cu}_{i_0+1}|_C$.
By the definition of $j_0$, $C\bigcap P^*_j=\phi$ for $j<i_0-1$,
so from proposition \ref{2.5}, $E^{cs}_{i_0-1}|_C$ is hyperbolic,
we denote it $E^s_{i_0-1}(C)$.                  \qed

\section{Proof of the main lemma}
\Pf At first, we can suppose $j_0\neq 0$, since if $j_0=0$, by
lemma \ref{4.3}, $C$ is a homoclinic class containing index 1
periodic points and $C$ is an index 0 fundamental limit, then we
proved the main lemma.

From lemma \ref{2.2}, there always exists a minimal index $j_0$
fundamental limit in $C$, we denote one of them $\Lambda$, by
lemma \ref{4.1}, we can suppose $\Lambda$ is a non-trivial minimal
set with a partial hyperbolic splitting $E^s_{j_0}\oplus
E^c_1\oplus E^u_{j_0+2}|_\Lambda$.

At first, let's prove $C$ contains an index $j_0$ or $j_0+1$
periodic point.

Now we divide the proof into two cases: $E^c_1(\Lambda)$ has an
$f$ orientation or not.

\noindent{\bf Case A:} $E^c_1(\Lambda)$ has an $f$ orientation.

At first, like in $\S$4, choose $V_0$ a small neighborhood of
$\Lambda$ such that $\Lambda_0=\bigcap f^i(\overline{V_0})$ will
have also a partial hyperbolic splitting $E^s_{j_0}\oplus E^c_1
\oplus E^u_{j_0+2}|_{\Lambda_0}$ and $E^c_1|_{\Lambda_0}$ also has
an $f$ orientation, and when $V_0$ is small enough, we can always
suppose the splitting can be extended to $V_0$ (of course, it's
not invariant any more). Choose $a_0$ small enough such that for
$x\in V_0$, we have $C^s_{a_0}(x)\bigcap C^{cu}_{a_0}(x)=\phi$,
$C^{cs}_{a_0}(x)\bigcap C^u_{a_0}(x)=\phi$, and
$C^i_{a_0}(x)\bigcap C^j_{a_0}(x)=\phi$ for ${(i\neq j\in
s,c,u)}$.

Choose another small neighborhood $V_1$ of $\Lambda$ satisfying
$\overline{V_1}\subset V_0$, let $\Lambda_1=\bigcap
\limits_{i=-\infty}^{\infty}f^i(\overline{V_1})$, then by lemma
\ref{4.5}, $\Lambda_1$ has a family of cental curves with uniform
size and locally invariant. Since $E^c_1(\Lambda_1)$ has an $f$
orientation, we choose one orientation and call the direction
right, at first let's consider the right central curves. By lemma
\ref{4.6}, we can suppose the family of right cental curves have
1-step contracting or expanding property (if it's not, by (a) of
lemma \ref{4.6} and lemma \ref{4.7}, $C$ is a homoclinic class
containing index $j_0$ or $j_0+1$ periodic point and $C$ is an
index $j_0$ fundamental limit.) that means for every point $x\in
\Lambda_1$, there exists a smooth central curve $\gamma^+_x\subset
l^+_{\delta_1}(x)$ on the right of $x$ such that
\begin{itemize}
\item $\gamma^+_x$ continuously depends on $x\in \Lambda_1$,

\item $f(\overline{\gamma^+_x})\subset \gamma^+_{f(x)}$ for all
$x\in \Lambda_1$ or $f^{-1}(\overline{\gamma^+_x})\subset
\gamma^+_{f^{-1}(x)}$ for all $x\in \Lambda_1$,

\item there exist $\varepsilon_0$ such that
$length(\gamma^+_x)>\varepsilon_0$.
\end{itemize}

Since $\Lambda$ is minimal, by 4) of proposition \ref{3.1}, there
exists a family of periodic points $\{p_n\}_{n=1}^\infty$ such
that $\lim \limits_{n\rightarrow\infty}Orb(p_n)\longrightarrow
\Lambda$, we can suppose $Orb(p_n)\subset V_1$ for $n\geq 1$, that
means $Orb(p_n)\subset \Lambda_1$, then on the right of $p_n$, we
have a central curve $\gamma^+_{p_n}$ such that either
$f(\overline{\gamma^+_{p_n}})\subset \gamma^+_{f(p_n)}$ (if the
right central curves of $\Lambda_1$ is 1-step contracting) or
$f^{-1}(\overline{\gamma^+_{p_n}})\subset \gamma^+_{f^{-1}(p_n)}$
(when the right central curves of $\Lambda_1$ is 1-step
expanding). For simplicity, we denote the central curve
$\gamma^+_{p_n}$ by $\gamma^+_n$, so we have
$f^{\pi(p_n)}(\overline{\gamma^+_n})\subset \gamma^+_n$ or
$f^{-\pi(p_n)}(\overline{\gamma^+_n})\subset \gamma^+_n$. Let
$\Gamma^+_n=\bigcap
\limits_{i=-\infty}^{\infty}f^{i\pi(p_n)}(\gamma^+_n)$, then
$\Gamma^+_n$ is a periodic segment with period $\pi(p_n)$, let
$q^+_n$ be one of the extreme point of $\Gamma^+_n$ different with
$p_n$ when $\Gamma^+_n$ is not trivial and $q^+_n=p_n$ when
$\Gamma^+_n$ is trivial, let $h^+_n=\gamma^+_n\setminus
\Gamma^+_n$, then there exists $\varepsilon_1$ doesn't depend on
$n$ such that $length(h^+_n)>\varepsilon_1$. It's easy to know
that when the right central curves are 1-step contracting,
$h^+_n\subset W^s(q^+_n)$ and if the right central curves are
1-step expanding, we have $h^+_n\subset W^u(q^+_n)$. With the same
argument on the left central curves, we can get $\gamma^-_n$,
$q^-_n$, $\Gamma^-_n$, $h^-_n$ also.

\begin{rk}\label{6.1}
 $\Gamma_n=\Gamma_n^+\bigcup \Gamma_n^-$ is a periodic central
 segment with period $\pi(p_n)$ and $f|_{\Gamma_n}$ is Kupka-Smale
diffeomorphism, that means $\Gamma_n$ just has finite fixed points
and they are sinks or sources.
\end{rk}

Now considering the contracting or expanding properties of the two
half parts of cental curves, we divide the proof into three
subcases:
\begin{itemize}
\item[A.1] Two sides of central curves are 1-step contracting.

\item[A.2] Right central curves are 1-step expanding and the left
central curves are 1-step contracting.

\item[A.3]Two sides of central curves are 1-step expanding.
\end{itemize}

\noindent{\bf Subcase A.1:} Two sides of central curves are 1-step
contracting.

In this subcase, we can show there exists periodic point $p\in C$
with index $j_0$ or $j_0+1$ and $Orb(p)\subset V_0$.\\

 We have known that $\gamma^+_{f^i(p_n)}\subset
l^+_{\delta_1}(f^i(p_n))\subset B_{\delta_1}(f^i(p_n))\subset
V_0$, and by 1-step contracting property , we have
$f^i(\gamma_n)\subset \gamma_{f^i(p_n)}$ for $i\geq 0$, so
$\gamma_n\subset \Lambda^+$, that means any $x\in \gamma_n$ has
uniform size of $\delta_1$ strong stable manifold
$W^{ss}_{\delta_1}(x)$. Since $\gamma_n$ is a positive central
segment, by the property of normally hyperbolic manifold and
$length(\gamma_n)>\varepsilon_0$ for all $n$, there exists
$\delta$ such that $W^s(\gamma_n)=\bigcup \limits_{x\in
\gamma_n}W^{ss}(x)$ is a cs disk with uniform size $\delta$, it's
easy to know $W^s(\gamma_n)=\bigcup \limits_{p\in
Per(\gamma_n)}W^s(p)$.

Let's suppose $\lim \limits_{n\rightarrow\infty}p_n=x_0\in
\Lambda$, then there exists $n$ big enough, such that
$W^{uu}(x_0)\h W^{cs}(\gamma_n)\neq \phi$, suppose $a\in
W^{uu}(x_0)\h W^{cs}(\gamma_n)$, by lemma \ref{Lya} we have $a\in
W^{uu}(x_0)\subset C$, it's easy to know $a\in W^s(p)$ for some
$p\in Per(\gamma_n)$, so $p\in \omega(a)\subset C$, recall that
all the central curves are in $V_0$, so $Orb(p)\subset V_0$ and
$p$ has index $j_0$ or
$j_0+1$.\\

\noindent{\bf Subcase A.2:} Right central curves are 1-step
expanding and the left central curves are 1-step contracting.\\

In this subcase, we can show that {\begin{itemize} \item[(a)]
either there exists periodic point $p\in C$ with index $j_0$ or
$j_0+1$ and $Orb(p)\subset V_0$

\item[(b)] or there exists periodic point $p\in C$ with index
$j_0$.
\end{itemize}}

From now we suppose that (a) is false, we claim that we can always
suppose $\lim \limits_{n\rightarrow
\infty}length(\Gamma^+_n)\longrightarrow 0$.

\noindent{\bf Proof of the claim:} Suppose there exist
$\delta^\prime$ and $\{\Gamma^+_{n_i}\}_{i=0}^{\infty}$ such that
$length(\Gamma^+_{n_i})>\delta^\prime$ for $i\geq 0$, then
$\gamma^-_{n_i}\bigcup \Gamma^+_{n_i}$ has uniform size and is a
positive central curves and $f^j(\gamma^-_{n_i}\bigcup
\Gamma^+_{n_i})\subset V$ for any $j\geq 0$, so like the argument
in Case A.1, $W^s(\gamma^-_{n_i}\bigcup \Gamma^+_{n_i})$ with
center $p_{n_i}$ has uniform size and when $i$ big enough, we have
$W^{uu}(x_0)\h W^s(\gamma^-_{n_i}\bigcup \Gamma^+_{n_i})\neq \phi$
, then $C$ contains an index $j_0$ or $j_0+1$ periodic point $p$
with $Orb(p)\subset V_0$, that's a contradiction with our
assumption that a) is false. \qed

Recall that central curves are a family of $C^1$ curves continuous
depend on $x\in \Lambda_1$, so we know $\lim \limits_{n\rightarrow
\infty}\gamma^+_n\longrightarrow \gamma_{x_0}^+$, with
$length(\Gamma^+_n)\longrightarrow 0$ we can know $\lim
\limits_{n\rightarrow\infty}h^+_n\longrightarrow \gamma^+_{x_0}$.

Since the right central curves are 1-step expanding we can know
$f^{-i}(h^+_n)\subset V_1$ for all $i\geq 0$, so
$f^{-i}(\gamma^+_{x_0})\subset V_0$ for all $i\geq 0$, that means
$\gamma^+_{x_0}$ is a negative central segment. With
$length(\gamma^+_{x_0})>\varepsilon_0$, that means
$W^u(\gamma^+_{x_0})=\bigcup \limits_{x\in
\gamma^+_{x_0}}W^{uu}_{\delta_1}(x)$ is a cu disk.

We claim that $\gamma^+_{x_0}\subset C$.\\
\noindent{\bf Proof of the claim:} Since $C$ is Lyapunov stable,
that means that there exists a family of open neighborhood
$\{U_n\}_{n=1}^\infty$ of $C$ such that
\begin{itemize}
\item[1)] $\overline{U_{n+1}}\subset U_n$

\item[2)] $f(\overline{U_n})\subset U_n$

\item[3)] $\bigcap \limits_n U_n =C$.
\end{itemize}

By the property of $\lim
\limits_{n\rightarrow\infty}Orb(p_n)=\Lambda$, we can suppose
$Orb(p_n)\subset U_n$ always, since $\lim \limits_{n\rightarrow
\infty}\Gamma^+_n\longrightarrow 0$, we can suppose $q^+_n\in U_n$
also. By the property of 2) above, we can know that
$W^u(q_n^+)\subset U_n$, since we've known that $h^+_n\subset
W^u(q_n^+)$, so $h^+_n\subset U_n$, then $\gamma^+_{x_0}=\lim
\limits_{n\rightarrow \infty}\gamma^+_n\subset
\bigcap\limits_{n\geq 1}U_n=C$.  \qed

\begin{rk}\label{central}
By the above argument, in fact we can know that for any $x\in
\Lambda$, $\gamma^+_{x}\subset C$. Then from
$f^{-i}(\gamma^+_{x})\subset \gamma^+_{f^{-i}(x)}$ for $i\in
\mathbb{N}$ we know $\gamma^+_{x}$ is a negative central segment,
so by e) of lemma \ref{4.4}, $\gamma^+_{x}$ has unstable manifold
$W^u(\gamma^+_{x})$, and by lemma \ref{Lya},
$W^u(\gamma^+_{x})\subset C$.
\end{rk}

Choose $y\in \gamma^+_{x_0}\setminus x_0$, then $y\in C$ also. Now
we claim that we can always suppose $y\in W^u(\Lambda)$.

\noindent {\bf Proof of the claim}: At first let's note that $y\in
\gamma^+_{x_0}\setminus x_0$ for $x_0\in \Lambda$ and
$f^{-i}(\gamma^+_{x_0})\subset \gamma^+_{f^{-i}(x_0)}$ for $i\in
\mathbb{N}$ because the right central model is 1-step expanding.
So $f^{-i}(y)\in \gamma^+_{f^{-i}(x_0)}$ for $i\in \mathbb{N}$,
that also means $\alpha(y) \subset V_1$ and $\alpha(y)$ has
partial hyperbolic splitting $E^s_{j_0}\oplus E^c_1 \oplus
E^u_{j_0+2}$. Now we just need show the length of the curve in
$\gamma^+_{f^{-i}(x_0)}$ which connecting $f^{-i}(y)$ and
$f^{-i}(x_0)$ converges to 0.

Now we suppose the length doesn't converge to 0, that means there
exists $i_n\longrightarrow \infty$ such that the length of the
curve in $\gamma^+_{f^{-i_n}(x_0)}$ which connecting $f^{-i_n}(y)$
and $f^{-i_n}(x_0)$ doesn't converge to 0, suppose $\lim
\limits_{n\rightarrow \infty}f^{-i_n}(x_0)\longrightarrow x_1$,
then $\lim \limits_{n\rightarrow \infty}f^{-i_n}(y)\longrightarrow
y_1\in \gamma^+_{(x_1)}\setminus x_1$. Since $y_1\in \alpha(y)$
and $\alpha(y)$ is a chain recurrent set in $V_1$, by generic
property 4) of proposition \ref{3.1}, there exists a family of
periodic orbits $\{Orb(p_n)\}\subset V_1$ such that
$p_n\longrightarrow y_1$, it's easy to know that $Orb(p_n)$ has
index $j_0$ or $j_0+1$ and $Orb(p_n)$ has uniform size of strong
stable manifold $W^{ss}_{\delta}(p_n)$, by e) of lemma \ref{4.5},
we know $W^{ss}_{\delta}(p_n)\bigcap W^u(\gamma^+_{x_1})\ni a\neq
\phi$.

Remark \ref{central} has shown that $a\in C$, so $Orb(p_n)\subset
\omega(a)\subset C$, recall that $\{Orb(p_n)\}\subset V_1$, then
we proved (a), it's a contradiction with the assumption that a) is
false. \qed

By the technique lemma, there exists a family of periodic points
$\{q_n\}$ such that $\lim \limits_{n\rightarrow \infty}q_n=y$ and
$y\bigcup\Lambda\subset \lim \limits_{n\rightarrow
\infty}Orb(q_n)$. By the corollary \ref{3.3}, we can suppose
$\{q_n\}$ all have index $j_0$ or index $j_0+1$. Denote
$C_0=\lim\limits_{n\rightarrow \infty}Orb(p_n)$, then $C_0\subset
P_{j_0}^*\bigcap P_{j_0+1}^*$, hence $C_0$ has a dominated
splitting $E^s_{j_0}\oplus E^c_1\oplus E^{cu}_{j_0+2}$, then by e)
of lemma \ref{4.5}, $q_n$ has uniform size of strong stable
manifold $W^{ss}_{\delta_1}(q_n)$ tangent at $q_n$ with
$E^s_{j_0}(q_n)$, and when $n$ big enough, we have
$W^{ss}_{\delta_1}(q_n)\h W^{u}(\gamma^+_{x_0})\neq \phi$. Let
$a\in W^{ss}_{\delta_1}(q_n)\h W^{u}(\gamma^+_{x_0})$, then $a\in
W^u(\gamma^+_{x_0})=\bigcup \limits_{x\in
\gamma^+_{x_0}}W^{uu}_{\delta_1}(x)\subset C$ and $q_n\in
\omega(a)\subset C$, so $C$ contains an index $j_0$ and index
$j_0+1$ periodic
point.\\

 \noindent{\bf Subcase A.3:} Two sides of central curves
are 1-step expanding.\\

In this subcase, we can show that
\begin{itemize}
\item[(a)] either there exists periodic point $p\in C$ with index
$j_0$ or $j_0+1$ and $Orb(p)\subset V_0$

\item[(b)] or there exists periodic point $p\in C$ with index
$j_0$.
\end{itemize}

We claim that if (a) is false, we can suppose there exists a
subsequence $\{n_j\}$ such that $\lim \limits_{j\rightarrow\infty}
\Gamma^+_{n_j}\longrightarrow 0$ or $\lim
\limits_{j\rightarrow\infty} \Gamma^-_{n_j}\longrightarrow 0$.

\noindent{\bf Proof of the claim:} Suppose it's wrong, then there
exists $\varepsilon>0$ such that $length(\Gamma^+_n)>\varepsilon$
and $length(\Gamma^-_n)>\varepsilon$, then $\Gamma_n$ is a central
segment with uniform size , with the same argument in case A.1, we
can show $C$ contains index $j_0$ or index $j_0+1$ periodic point,
and its orbit is contained in $V_0$.\qed

Now change by a subsequence, we can suppose $\lim
\limits_{n\rightarrow\infty} \Gamma^+_{n}\longrightarrow 0$, now
the rest argument is the same with case A 2.\\

\noindent{\bf Case B: $E^c_1(\Lambda)$ has no $f$ orientation:}\\

In this case, we can locally define orientation, and in this case
locally the two sides of central curves are either 1-step
expanding or 1-step contracting, the rest argument is almost the
same with Case A.1 and Case A.3.\\

Now let's prove that $C$ is an index $j_0$ fundamental limit, here
we choose a family of neighborhoods $\{V_n\}$ of $\Lambda$ such
that $V_{n+1}\subset V_n$ and $\bigcap V_n=\Lambda$, then by above
argument, we can show that
\begin{itemize}
\item[(a)] either $C$ contains index $j_0$ periodic point,

\item[(b)] or $C$ contains periodic point $p_n\in C$ with index
$j_0+1$ and $Orb(p_n)\subset V_n$.
\end{itemize}
In the case (a), of course $C$ is an index $j_0$ fundamental
limit; in the case (b) we just need the following lemma given in
\cite{1Y2}:

\begin{lem}\label{fund}
Suppose $f\in R$, $C$ is a non-trivial chain recurrent class of
$f$, and $\Lambda\varsubsetneq C$ is a minimal set with partial
hyperbolic splitting $E^s_{j_0}\oplus E^c_1\oplus E^u_{j_0+2}$
where $dim(E^c_1(\Lambda))=1$ and $E^c_1(\Lambda)$ is not
hyperbolic, if there exists a family of periodic points $\{p_n\}$
in $C$ satisfying $\lim \limits_{n\rightarrow
\infty}Orb(p_n)=\Lambda$, then $C$ is index $j_0$ and $j_0+1$
fundamental limit.
\end{lem}

\begin{rk}
The proof of the above lemma is divided into two cases:
\begin{itemize}
\item[(A)] there exists $\delta>0$ such that for any $p_n$, we
have $Df^{\pi(p_n)}|_{E^c_1(p_n)}<e^{-\delta\pi(p_n)}$,

\item[(B)] for any $\frac{1}{m}$, there exists $p_{n_m}$ such that
$Df^{\pi(p_{n_m})}|_{E^c_1(p_n)}>e^{-\frac{1}{m}\pi(p_{n_m})}$.
\end{itemize}

In the first case we use weakly selecting lemma, and in case (B)
we use lemma 4.25 of \cite{1Y2} which basically is a transition
property.
\end{rk}             \qed

\section{Proof of technique lemma}

The proof of the technique lemma depends on generic assumption
heavily, with many generic assumptions, we can find some segment
of orbit with 'good' position, then after using connecting lemma
and another generic property, we can get the periodic points which
we need.

In $\S$7.1, we'll introduce some new $C^1$ generic properties in
order to define the generic set given in technique lemma. In
$\S$7.2, we'll recall the proof of connecting lemma, especially
about the 'cutting tool', because we need an important fact which
just appears in the proof of connecting lemma. In $\S$7.3, we'll
prove the technique lemma.

\subsection{Some new $C^1$ generic properties}.
Suppose $\{U_\alpha\}_{\alpha\in \mathcal A}$ is a topological
basis of $M$ satisfying for any $\varepsilon>0$, there exists a
subsequence $\{U_{\alpha_i}\}_{i=1}^\infty$ such that
$diam(U_{\alpha_i})<\varepsilon$ and $\bigcup
\limits_{i}(U_{\alpha_i})$ is a cover of $M$. Fix this topological
basis, we'll get some new $C^1$ generic properties.

At first, let's recall some definitions, suppose $K$ is a compact
set of $M$, $f\in C^1(M)$ has been given, $x,y\in K$,
$x\underset{K}{\dashv}y$ means that for any $\varepsilon>0$, there
exists an $\varepsilon$-pseudo orbit in $K$ beginning from $x$ and
ending at $y$. If $K=M$, we just denote $x\dashv y$.

The following result has been proved in [Cr2]:

\begin{lem}\label{7.1}
There exists a generic subset $R^*_{1,0}$ such that any $f\in
R^*_{1,0}$ will satisfy the following property: suppose $K$ is a
compact set, $W$ is any neighborhood of $K$, $x_1,x_2\in K$
satisfy $x_1\underset {K}{\dashv}x_2$, $U_1,U_2\subset W$ are
neighborhoods of $x_1,x_2$ respectively, then there exists a
segment of orbit of $f$ in $W$ beginning from $U_1$ and ending in
$U_2$. More precisely, there exists $a\in U_1$ and $i_1>0$ such
that $f^{i_1}(a)\in U_2$ and $f^i(a)\in W$ for $0\leq i \leq i_1$.
\end{lem}

\begin{lem}\label{7.2}
There exists a generic subset $R^*_{1,1}$ such that any $f\in
R^*_{1,1}$ will satisfy the following property: suppose $\Lambda$
is a invariant compact subset of $f$, $y\notin \Lambda$, $0<s<1$,
$\{\Phi_i\}_{i=1}^K\subset\{U_\alpha\}_{\alpha\in \mathcal A}$ is
an open cover for $\Lambda$ and $O\in \{U_\alpha\}_{\alpha\in
\mathcal A}$ is a small neighborhood of $y$, if there exist
$g_n\overset{c^1}{\longrightarrow}f$ and $g_n$ has periodic point
$p_n$ satisfying $\frac{\#\{Orb_{g_n}(p_n)\bigcap (\bigcup
\limits_{i=1}^K\Phi_i)\}}{\pi_{g_n}(p_n)}>s$ and
$Orb_{g_n}(p_n)\bigcap O\neq \phi$, then $f$ itself has a periodic
point $p$ satisfying $\frac{\#\{Orb(p)\bigcap (\bigcup
\limits_{i=1}^K\Phi_i\})}{\pi(p)}>s$ and $Orb(p_n)\bigcap O\neq
\phi$.
\end{lem}

\Pf: Consider the set
$\{(\Phi_{\beta_1},\cdots,\Phi_{\beta_{N(\beta)}};O_{\beta})\}
_{\beta \in {\mathcal B}_0}$ where $\Phi_{\beta_i},O_\beta\in
\{U_\alpha\}_{\alpha\in \mathcal A}$, it's easy to know ${\mathcal
B}_0$ is countable.

For any $\beta \in {\mathcal B}_0$, denote
\begin{itemize}
\item $H_\beta=\{f|\ f\in C(M), f$ has a $C^1$ neighborhood
$\mathcal U$ such that for any $g\in \mathcal U$, $g$ has a
periodic orbit $p_g$ satisfying $\frac{\#\{Orb_{g}(p_g)\bigcap
(\bigcup \limits_{i=1}^K\Phi_i)\}}{\pi_{g}(p_g)}>s$ and
$Orb_{g}(p_g)\bigcap O_\beta\neq \phi\}$,

\item $N_\beta=\{f|\ f\in C^1(M), f$ has a $C^1$ neighborhood
$\mathcal U$ such that for any $g\in \mathcal U$, $g$ has no any
periodic orbit $p_g$ satisfying $\frac{\#\{Orb_{g}(p_g)\bigcap
(\bigcup \limits_{i=1}^K\Phi_i)\}}{\pi_{g}(p_g)}>s$ and
$Orb_{g}(p_g)\bigcap O_\beta\neq \phi\}$.
\end{itemize}

It's easy to know $H_\beta \bigcup N_\beta$ is open and dense in
$C^1(M)$. Let $R^*_{1,0}=\bigcap \limits_{\beta \in {\mathcal
B}_0}(H_\beta \bigcup N_\beta)$, we'll show $R^*_{1,0}$ satisfies
the property we need.

For any $f \in R^*_{1,0}$ and any $\beta^*\in {\mathcal B}_0$,
suppose there exists a family of $C^1$ diffeomorphisms
$\{g_n\}^\infty_{n=1}$ such that $\lim
\limits_{n\rightarrow\infty}g_n=f$ and any $g_n$ has a periodic
orbit $p_n$ satisfying $\frac{\#\{Orb_{g_n}(p_n)\bigcap (\bigcup
\limits_{i=1}^K\Phi_i)\}}{\pi_{g_n}(p_n)}>s$ and
$Orb_{g_n}(p_n)\bigcap O_\beta\neq \phi$, then $f\notin
N_{\beta^*}$. That means $f \in H_{\beta^*}$, so we proved this
lemma. \qed

With the same argument like above, we can get the following result
:
\begin{lem}\label{avoid}
There exists a generic subset $R_{1,2}^*$ such that any $f\in
R_{1,2}^*$ will satisfy the following property: for finite number
of open set $\{\Phi_i\}_{i=1}^{N}\subset\{U_\alpha\}_{\alpha\in
\cal A}$ and $U_0,U_1,U_2,U_3\in\{U_\alpha\}_{\alpha\in\cal A}$
such that $U_0,U_1,U_2,U_3\subset \bigcup
\limits_{i=1}^{N}\Phi_i$, if there exist $a_n\in U_0$,
$g_n\overset{C^1}{\longrightarrow}f$ and $0<i_{1,n}<i_{2,n}$ such
that $g_n^i(a_n)\in \bigcup \limits_{i=1}^{N}\Phi_i$ for $0\leq i
\leq i_{2,n}$, $g_n^{i_{1,n}}(a_n)\in U_1$, $g_n^{i_{2,n}}(a_n)\in
U_3$ and $g_n^{i}(a_n)\notin \overline{U_2}$ for $0\leq i \leq
i_{1,n}$, then there exist $a\in U_0$ and $0<i_1<i_2$ such that
$f^{i}(a)\in \bigcup \limits_{i=1}^{N}\Phi_i$ for $0\leq i \leq
i_2$, $f^{i_1}(a)\in U_1$, $f^{i_2}(a)\in U_3$ and $f^i(a_n)\notin
\overline{U_2}$ for $0\leq i\leq i_1$.
\end{lem}

Now let $R_0^\prime=R_0\bigcap R_{1,0}^*\bigcap R_{1,1}^*\bigcap
R_{1,2}^*$, and in $\S$7.3 we'll show the set will satisfy the
technique lemma.

\subsection{Introduction of connecting lemma}

Connecting lemma was proved by Hayashi \cite{1H} at first, and
then was extended to the conservative setting by Xia, Wen
\cite{1XW}. the following statement of connecting lemma was given
by Lan Wen as an uniform version of connecting lemma.

\begin{lem} \label{7.3} (connecting lemma \cite{1W2}) For any $C^1$ neighborhood
$\mathcal U$ of $f$, there exist $\rho>1$, a positive integer $L$
and $\delta_0>0$ such that for any $z$ and $\delta<\delta_0$
satisfying
$\overline{f^i(B_\delta(z))}\bigcap\overline{f^j(B_\delta(z))}=\phi$
for $0\leq i\neq j\leq L$, then for any two points $p$ and $q$
outside the cube $\Delta=\bigcup^L_{i=1}f^i(B_\delta(z))$, if the
positive $f$-orbit of $p$ hits the ball $B_{\delta/\rho}(z)$ after
$p$ and if the negative $f$-orbit of $q$ hits the small ball
$B_{\delta/\rho}(z)$, then there is $g\in \mathcal U$ such that
$g=f$ off $\Delta$ and $q$ is on the positive $g$-orbit of $p$.
\end{lem}

\begin{rk}\label{connecting}
Suppose we have another point $z_1\in M$ satisfying
$\Delta_1\bigcap \Delta=\phi$ where
$\Delta_1=\bigcup^L_{i=1}f^i(B_\delta(z_1))$, then if we use twice
connecting lemma in $\Delta$ and $\Delta_1$, we can still get a
diffeomorphism $g$ in $\cal U$.
\end{rk}

Now we'll show the idea of the proof of connecting lemma, because
we need some special property which just appears in the proof.

In the proof, the main idea is Hayashi's 'cutting' tool, by it we
can cut some orbits from $p$'s original $f$-orbit and $q$'s
original $f$-orbit, and then connect the rest part in $\Delta$.
More precisely description is following. Suppose $f^{s_m}(p)\in
B_{\delta/\rho}(z)$ and there exists $0<s_1<s_2<\cdots<s_m$ such
that $f^{s_i}\in B_\delta(z)$ for $1\leq i \leq m$ and
$f^s(p)\notin B_\delta(z)$ for $s\in
\{0,1,\cdots,s_m\}\setminus\{s_1,s_2,\cdots,s_m\}$. For $q$, there
exists $0<t_1<t_2<\cdots<t_n$ such that $f^{-t_i}(q)\in
B_\delta(z)$ for $1\leq i \leq n$, $f^{t_n}(q)\in
B_{\delta/\rho}(z)$ and $f^{-t}(q)\notin B_\delta(z)$ for
$t\in\{0,1,\cdots,t_n\}\setminus\{t_1,t_2,\cdots,t_n\}$. By some
rule, we can cut some $f$-orbits in $p$'s orbit like
$\{f^{s_i+1}(p),f^{s_i+2}(p),\cdots, f^{s_j}(p)\}_{j>i}$ and cut
some $f$-orbits in $q$'s orbit like
$\{f^{-t_i}(q),\cdots,f^{-t_j+2}(q),f^{-t_j+1}(q)\}_{j>i}$, the
rest segment is like:
       $$P^\prime=(p,f(p),\cdots,f^{s_{i_1}}(p);f^{s_{i_2}+1},
       \cdots,f^{s_{i_3}}(p);\cdots;f^{s_{i_{(k(p)-1)}}+1}(p),\cdots,f^{s_{i_{k(p)}}}(p)),$$
       $$Q^\prime=(
       f^{-t_{j_{k(q)}}+1}(q),\cdots,f^{-t_{j_{k(q)-1}}}(q);\cdots;f^{-t_{j_3}+1},\cdots,
       f^{-t_{j_2}}(q);f^{-t_{j_1}+1}(q)\cdots,f^{-1}(q),q,).$$
Denote $X=P^\prime \bigcup Q^\prime$, and $\pi(X)$ is the length
of $X$, it's easy to know $X$ is a $2\delta$-pseudo orbits. Then
we can do several perturbations called 'push' in $\Delta$ and get
a diffeomorphism $g$ such that $q$ is on the positive $g$-orbit of
$p$, in fact, we have $g^{\pi(X)}(p)=q$. It's because after the
push, we can connect $f^{s_{i_1}}(p)$ and $f^{s_{i_2}+L}(p)$,
$\cdots$; $f^{s_{i_{k(p)-2}}}(p)$ and $f^{s_{i_{k(p)-1}}+L}(p)$;
$f^{s_{i_{k(p)}}}(p)$ and $f^{-t_{j_{k(q)}}+L}(q)$;
$f^{-t_{j_{k(q)-1}}}(q)$ and $f^{-t_{j_{k(q)-2}}+L}(q)$; $\cdots$;
$f^{-t_{j_2}}(q)$ and $f^{-t_{j_1}+L}(q)$ by $L$ times pushes in
$\Delta$, we don't cut orbits anymore, and it's important to note
that the supports of different pushes don't intersect with each
other, so we don't change the length of $X$, we just push the
points of $X$ in $\Delta$ and get a connected orbit. By the above
argument, it's easy to know $g|_{M\setminus \Delta}=f|_{M\setminus
\Delta}$ and $g(\Delta)=f(\Delta)$.

\begin{rk}\label{7.4} In the above argument, suppose there exists an open set $V$ such that
$f^i(p)\in V$ for $0\leq i \leq s_m$ and $\Delta \subset V$, then
after cutting and pushes, we can know
$\{p,g(p),\cdots,g^{\pi(P^\prime)}(p)\}\subset V$. What's more, we
can show that
$\#\{\{g^i(p)\}_{i=0}^{\pi(P^\prime)+\pi(Q^\prime)}\bigcap
(V)^c\}<t_n$.
\end{rk}

\subsection{Proof of technique lemma}

\Pf: Here we just prove the technique lemma for $y\in C\bigcap
W^s(\Lambda)\setminus C$, the proof for the other case is similar.

Fix $V_0$ a small neighborhood of $\Lambda$, $U_0\in
\{U_\alpha\}_{\alpha\in \cal A}$ a small neighborhood of $y$ such
that $\overline{V_0}\bigcap \overline{U_0}=\phi$ and
$\overline{U_0}\subset O$, $\overline{V_0}\subset V$. Let
$\{\Phi_i\}^N_{i=1}\subset\{U_\alpha\}_{\alpha\in \cal A}$ be an
open cover of $\Lambda$ such that $\overline{\bigcup
\limits_{i=1}^{N}\Phi_i}\subset V_0$, choose $V_1$ another small
compact neighborhood of $\Lambda$ such that $V_1\subset \bigcup
\limits_{i=1}^{N}\Phi_i$, and choose
$O(y)\in\{U_\alpha\}_{\alpha\in \cal A}$ is a neighborhood of $y$
such that $\overline{O(y)}\subset U_0$.

Choose $x_0\in \omega(y)\subset \Lambda$, denote $z_0=f^{i_0}(y)$
is the last time the positive orbit of $y$ enters $V_1$, then we
have $z_0\underset{V_1}{\dashv}x_0$. It's easy to know that $z_0$
is not a periodic point.

\begin{figure}[h]
\begin{center}
\psfrag{U0}{$U_0$}\psfrag{V0}{$V_0$}\psfrag{V1}{$V_1$}\psfrag{y}{$y$}\psfrag{x0}{$x_0$}\psfrag{z0}{$z_0$}\psfrag{z1}{$z_1$}\psfrag{Lambda}{$\Lambda$}
\includegraphics[height=1.5in]{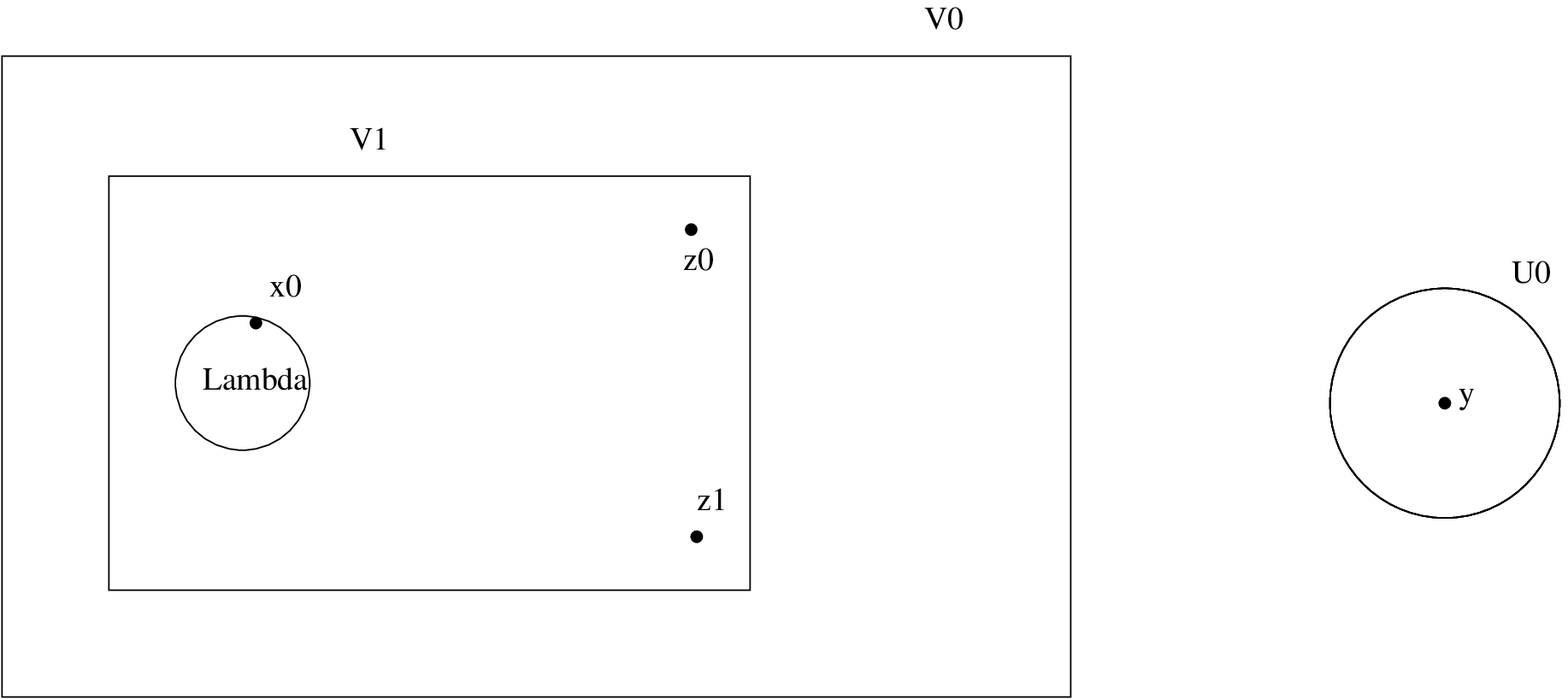}
\end{center}
\end{figure}

Now choose $\{\delta_n\}_{n=1}^{\infty}$ satisfying
$\delta_n\longrightarrow 0^+$, for every $\delta_n$, there exists
a $\delta_n$-pseudo orbit from $x_0$ to $y$, denote $z^-_n$ the
first time the pseudo orbit leaves $V_1$, suppose $\lim
\limits_{n\rightarrow \infty}f^{-1}(z_n^-)=z_1$, then
$Orb^-(z_1)\in \overline{V_1}$ and $x_0\underset{V_1}{\dashv}z_1$.
We can always suppose $z_1$ is not a periodic point, since if
$z_1$ is a periodic point, by $f$ is a Kupka-Smale diffeomorphism,
$z_1$ should be a hyperbolic periodic point, then there exists a
point $z_1^\prime\in W^s_{loc}(z)$ such that $Orb^-(z_1^\prime)\in
\overline{V_1}$ and $z_1^\prime\underset{V_1}{\dashv}x_0$, then we
can replace $z_1$ by $z_1^\prime$.

Before we enter the details of the proof, we'll show some ideas of
the proof. In the beginning we show that there exists a orbit
beginning from a neighborhood of $z_1$ to a neighborhood of $z_0$.
Then we show there exists another segment of orbit in
$\overline{V_1}$ beginning from a neighborhood of $z_0$ passing a
very small neighborhood of $x_0$ and ending in a neighborhood of
$z_1$, and most important, the orbit between the neighborhood of
$z_0$ and the neighborhood of $x_0$ will never pass $z_1$'s
neighborhood. Here we should note that until now we just use
generic property, and we don't do any perturbation yet. Now we'll
use connecting lemma twice to connect the above two orbits and get
a periodic orbit, more precisely, at first we use connecting lemma
at $z_1$'s neighborhood and then we use connecting lemma near
$z_0$'s neighborhood, and we can show after the perturbations, the
periodic orbit we get will spend a long time in $V_1$, then with
generic assumption again, we can know that $f$ itself has such
kind of periodic orbit.

At first, we need the following lemma which can help us obtain an
orbit with 'good' position:

\begin{lem}\label{7.6}
There exists $\delta_0>0$ such that for any
$\delta_1^*,\delta_2^*,\delta_3^*<\delta_4^*<\delta_0$, there
exist $a\in B_{\delta_1^*}(z_0)$ and $0<i_1<i_2$ such that
$f^{i_1}(a)\in B_{\delta_2^*}(x_0)$, $f^{i_2}(a)\in
B_{\delta_3^*}(z_1)$, $f^i(a)\in \bigcup \limits_{i=1}^{N}\Phi_i$
for $0\leq i \leq i_2$ and $f^i(a)\notin
\overline{B_{\delta_4^*}(z_1)}$ for $0\leq i \leq i_1$.
\end{lem}

\Pf: Since $Orb^-(z_1)\subset V_1$ and $f^{-i_0}(z_0)=y\notin
V_1$, we get $z_1\notin Orb^+(z_0)$, with the fact $z_1\notin
\omega(z_0)$, we can choose $Y^+(z_1)\in \{U_\alpha\}_{\alpha\in
\cal A}$ a small neighborhood of $z_1$ such that
$\overline{Y^+(z_1)}\bigcap Orb^+(z_0) =\phi$, $Y^+(z_1)\subset
\bigcup \limits_{i=1}^{N}\Phi_i$ and $\overline{Y^+(z_1)}\bigcap
\Lambda=\phi$. Choose $\delta_0>0$ small enough such that
\begin{itemize}
\item $B_{\delta_0}(z_1)\subset Y^+(z_1)$,

\item $B_{\delta_0}(z_0)\subset \bigcup \limits_{i=1}^{N}\Phi_i$,
$\overline{B_{\delta_0}(z_0)}\bigcap \Lambda=\phi$, and
$B_{\delta_0}(z_0)\bigcap Y(z_1)=\phi$,

\item $B_{\delta_0}(x_0)\subset \bigcup \limits_{i=1}^{N}\Phi_i$,
$B_{\delta_0}(x_0)\bigcap Y(z_1)=\phi$, $B_{\delta_0}(x_0)\bigcap
B_{\delta_0}(z_0)=\phi$.
\end{itemize}

Now suppose $\delta_1^*,\delta_2^*,\delta_3^*<\delta_4^*<\delta_0$
are fixed, we can choose $X(z_0)\in \{U_\alpha\}_{\alpha\in \cal
A}$ a small neighborhood of $z_0$ satisfying $X(z_0)\subset
B_{\delta_{1}^*}(z_0)$ and choose $Y^-(z_1)\in
\{U_\alpha\}_{\alpha\in\cal U}$ a small neighborhood of $z_1$ such
that $Y^-(z_1)\subset B_{\delta_3}^*(z_1)\subset Y^+(z_1)$. For
any small $\varepsilon_n>0$, by connecting lemma,
$B_{\varepsilon_n}(f)$ gives us parameters $L_n,\delta_n$ and
$\rho_n$, we choose $W_n^+,W_n^-\in \{U_\alpha\}_{\alpha\in \cal
A}$ neighborhoods of $x_0$ small enough such that

\begin{itemize}
\item $W^+_n,W^-_n\subset B_{\delta_{2}^*}(x_0)$,

\item there exists $0<\delta<\delta_n$ such that $W_n^-\subset
B_{\delta/\rho_n}(x_0)\subset B_{\delta}(x_0)\subset W_n^+\subset
B_{\delta_2^*}(x_0)$ and we have $f^i(W_n^+)\bigcap
f^j(W_n^-))=\phi$ for $0\leq i\neq j \leq L_n$.

\item denote $\Delta_n=\bigcup \limits_{i=0}^{L_n}W_n^+(x_0)$,
then $\Delta_n\subset \bigcup \limits_{i=1}^N\Phi_i$ and
$\Delta_n\bigcap X(z_0)=\phi$, $\Delta_n\bigcap Y^+(z_1)=\phi$.
\end{itemize}
Since $\Lambda$ is an invariant compact subset, $z_0,z_1\notin
\Lambda$ and $x_0\in \Lambda$ is not a periodic point, we can
always choose such kind of neighborhoods.

Since $x_0\in \omega(z_0)\subset\Lambda$, then there exists
$i_{1,n}^*$ such that $f^{i_{1,n}^*}(z_0)\in W_n^-$; because $x_0
\underset{V_1}{\dashv} z_1$, by lemma \ref{7.1}, there exist
$b_n\in Y^-(z_1)$ and $j_{n}$ such that $f^{-j_{n}}(b_n)\in W^-_n$
and $f^{-j}(b_n)\in \bigcup \limits_{i=1}^{N}\Phi_i$ for $0\leq j
\leq j_{n}$.

Recall $W_n^-\subset B_{\delta/\rho_n}(x_0)$, use connecting lemma
to connect $z_0$ and $b_n$ in $\Delta_n$, we can get a new
diffeomorphism $g_n$ and $i_{0,n},i_{1,n}$ such that
$g_n^{i_{0,n}}(a_n)\in W^+_n$, $g_n^{i_{1,n}}(a_n)=b_n\in
Y^-(z_1)$; since the original two orbits are both in $\bigcup
\limits_{i=1}^{N}\Phi_i$ and $\Delta_n\subset \bigcup
\limits_{i=1}^{N}\Phi_i$, we can know that $g_n^j(a_n)\in \bigcup
\limits_{i=1}^{N}\Phi_i$ for all $0\leq j \leq i_{1,n}$.

From $\{f^j(z_0)\}_{j=0}^{i_{1,n}^*}\subset
(\overline{Y^+(z_1)})^c$, by remark \ref{7.4}, we can choose
$i_{0,n}$ such that $(g_n)^j(a_n)\notin \overline{Y^+(z_1)}$ for
$0\leq j \leq i_{0,n}$.

Now fix $n_0\in \mathbb{N}$ and consider the neighborhood
$W^+_{n_0}$ of $x_0$, with generic property lemma \ref{avoid},
there exists $a\in X(z_0)$ and $0<i_1<i_2$ such that $f^i(a)\in
\bigcup \limits_{i=1}^{N}\Phi_i$ for $0\leq i \leq i_2$,
$f^{i_1}(a)\in W^+_{n_0}$, $f^{i_2}(a)\in Y^-(z_1)$ and
$f^{i}(a)\notin \overline{Y^+(z_1)}$ for $0\leq i \leq i_1$. With
the facts $X(z_0)\subset B_{\delta_{1}^*}(z_0)$, $W^+_n\subset
B_{\delta_{2}^*}(x_0)$ and $Y^-(z_1)\subset
B_{\delta_{3}^*}(z_1)\subset Y^+(z_1)$, we finish the proof. \qed

Now for any sequence $\varepsilon_n\longrightarrow0^+$ and
$s_n\longrightarrow1^-$, consider $B_{\varepsilon_n}(f)$ the
$\varepsilon_n$-neighborhood of $f$ in $C^1(M)$, by connecting
lemma $B_{\varepsilon_n}(f)$ gives us a family of parameters
$\rho_n\longrightarrow \infty$, $\delta_n\longrightarrow 0$ and
$L_n$. Then there exist $\delta_{0,n}\longrightarrow 0^+$ such
that
\begin{itemize}
\item[A1] $\delta_{0,n}<\delta_{n}$, $\delta_{0,n}<\delta_0$

\item[A2] $\delta_{0,n+1}<\delta_{0,n}/\rho_n$,

\item[A3] $f^i(B_{\delta_{0,n}}(z_0))\bigcap
f^j(B_{\delta_{0,n}}(z_0))=\phi$ for $0\leq i\neq j \leq L_n$ and
$\bigcup_{i=1}^{L_n}f^i(B_{\delta_{0,n}}(z_0))\subset
\bigcup_{i=1}^{N}\Phi_i$,

\item[A4] $\overline{B_{\delta_{0,n}}(z_0)}\bigcap \Lambda=\phi$
and $\overline{\bigcup
_{i=0}^{L_n}f^i(B_{\delta_{0,n}}(z_0))}\bigcap \bigcup
_{i=0}^{L_n}f^{-i}(z_1)=\phi$.
\end{itemize}
Since $z_0$ is not periodic point, $Orb^+(z_0)\subset
\overline{V_1}\subset \bigcup \limits_{i=1}^{N}\Phi_i$ and
$\omega(z_0)\subset \Lambda$, we can always choose the above
sequence $\{\delta_{0,n}\}$ for $z_0$. For $z_1$ we can also
choose a sequence $\{\delta_{1,n}\}$ such that
\begin{itemize}
\item[B1] $\delta_{1,n}<\delta_{n}$, $\delta_{1,n}<\delta_{0}$

\item[B2] $\delta_{1,n+1}<\delta_{1,n}/\rho_n$,

\item[B3] $f^{-i}(B_{\delta_{1,n}}(z_1))\bigcap
f^{-j}(B_{\delta_{1,n}}(z_1))=\phi$ for $0\leq i\neq j \leq L_n$
and $\bigcup_{i=0}^{L_n}f^{-i}(B_{\delta_{1,n}}(z_1))\subset
\bigcup_{i=1}^{N}\Phi_i$,

\item[B4] $\overline{B_{\delta_{1,n}}(z_1)}\bigcap \Lambda=\phi$
and $\overline{\bigcup
_{i=0}^{L_n}f^i(B_{\delta_{0,n}}(z_0))}\bigcap \overline{\bigcup
_{i=0}^{L_n}f^{-i}(B_{\delta_{1,n}}(z_1))}=\phi$.
\end{itemize}

\begin{figure}[h]
\begin{center}
\psfrag{V1}{$V_1$}\psfrag{Wn}{$B_{\delta_{2,n}}(x_0)$}\psfrag{Xn}{$B_{\delta_{0,n}}(z_0)$}
\psfrag{Yn}{$B_{\delta_{1,n}}(z_1)$}\psfrag{x0}{$x_0$}\psfrag{z0}{$z_0$}\psfrag{z1}{$z_1$}\psfrag{L}{$\Lambda$}
\includegraphics[height=1.5in]{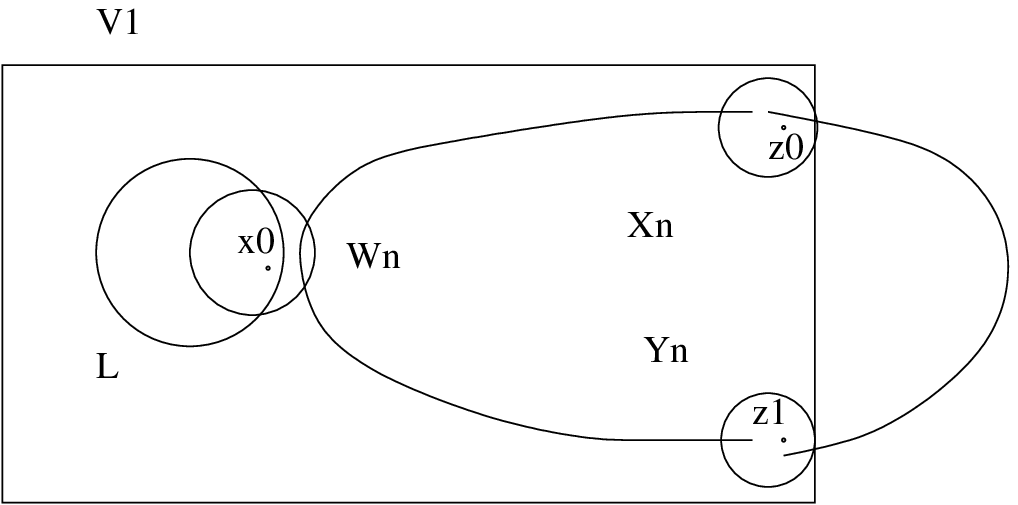}
\end{center}
\end{figure}

Then by lemma \ref{7.1}, there exists a family of points $\{a_n\}$
in $B_{\delta_{1,n}/\rho_n}(z_1)$ and $i_{0,n}$ such that
$f^{i_{0,n}}(a_n)\in B_{\delta_{0,n}/\rho_n}(z_0)$. We define
$\Delta_{0,n}=\bigcup
\limits_{i=1}^{L_n}f^i(B_{\delta_{0,n}}(z_0))$ and
$\Delta_{1,n}=\bigcup
\limits_{i=1}^{L_n}f^{-i}(B_{\delta_{1,n}}(z_1))$.

Now we'll choose a sequence of number
$\delta_{2,n}\longrightarrow0^+$ such that:
\begin{itemize}
\item[C1] $\delta_{2,n+1}<\delta_{2,n}$,
$\delta_{2,n}<\delta_{0}$,

\item[C2] $B_{\delta_{2,n}}(x_0)\subset V_1$,
$B_{\delta_{2,n}}(x_0)\bigcap \Delta_{0,n}=\phi$ and
$B_{\delta_{2,n}}(x_0)\bigcap \Delta_{1,n}=\phi$.

\item[C3] For any $j_0$ satisfying
$f^{j_0}(B_{\delta_{0,n}}(z_0))\bigcap B_{\delta_{2,n}}(x_0)\neq
\phi$, we have $\frac{i_{0,n}}{j_0}<1-s$.
\end{itemize}
Since $\Lambda$ is an invariant compact subset in $V_1$, we can
always choose such neighborhoods.

Now by lemma \ref{7.6}, for
$B_{\delta_{0,n}/\rho_n}(z_0),B_{\delta_{1,n}/\rho_n}(z_1)\subset
B_{\delta_{1,n}}(z_1)$ and $B_{\delta_{2,n}}(x_0)$ there exists an
orbit in $\bigcup \limits_{i=1}^{N}\Phi_i$ beginning in
$B_{\delta_{0,n}/\rho_n}(z_0)$ passing $B_{\delta_{2,n}}(x_0)$ and
ending in $B_{\delta_{1,n}/\rho_n}(z_1)$. More precisely, it means
that there exist $b_n\in B_{\delta_{1,n}/\rho_n}(z_0)$ and
$0<j_{0,n}^*\leq j_{1,n}^*<j_{2,n}^*$ such that:
\begin{itemize}
 \item[D1] $f^j(b_n)\in \bigcup \limits_{i=1}^N \Phi_i$ for $0\leq j \leq j_{2,n}^*$,

\item[D2] $f^{j_{0,n}^*}(b_n)\in B_{\delta_{0,n}}(z_0)$,
$f^{j_{1,n}^*}(b_n)\in B_{\delta_{2,n}}(x_0)$,
$f^{j_{2,n}^*}(b_n)\in B_{\delta_{1,n}/\rho_n}(z_1)$,

\item[D3] $f^j(b_n)\notin B_{\delta_{0,n}}(z_0)$ for $j_{0,n}^*< j
\leq j_{1,n}^*$, and $f^j(b_n)\notin
\overline{B_{\delta_{2,n}}(z_1)}$ for $0\leq j \leq j_{1,n}^*$.
\end{itemize}

\begin{rk}\label{7.7}
In fact, we can know that $\{f^j(b_n)\}_{j={0}}^{j_{1,n}^*}\bigcap
\Delta_{1,n}=\phi$ and
$\{f^j(b_n)\}_{j=j_{0,n}^*+L_n}^{j_{1,n}^*}\bigcap
\Delta_{0,n}=\phi$, so in the following proof, when we use
connecting lemma in $\Delta_{0,n},\Delta_{1,n}$ twice, we can get
a new diffeomorphism $g_n$ and a periodic orbit $Orb_{g_n}(p_n)$
of $g_n$ such that the segment
$\{f^j(b_n)\}_{j=j_{0,n}^*+L_n}^{j_{1,n}^*}\subset Orb_{g_n}(p_n)$
and $\#\{Orb_{g_n}(p_n)\bigcap (\bigcup
\limits_{i=1}^N\Phi_i)^c\}<i_{0,n}$, then by C3, we can know that
$\frac{\#\{Orb_{g_n}(p_n)\bigcap \bigcup
\limits_{i=1}^N\Phi_i\}}{\pi_{g_n}(p_n)}>1-\frac{i_{0,n}}{j_{1,n}^*-j_{0,n}^*}>s$.
\end{rk}

Now fix an $n$, let's consider the two points $f^{i_{0,n}}(a_n)$
and $b_n$, we know the positive $f$-orbit of $b_n$ hits
$B_{\delta_{1,n}/\rho_n}(z_1)$ after $b_n$ and the negative
$f$-orbit of $f^{i_{0,n}}(a_n)$ hits
$B_{\delta_{1,n}/\rho_n}(z_1)$ also, by connecting lemma, the fact
$\Delta_{1,n}\subset \bigcup \limits_{i=1}^{N}\Phi_i$, property D3
and remark \ref{7.4}, \ref{7.7}, there exists $g^*_n\in
B_{\varepsilon_n}(f)$ such that $g^*_n\equiv f$ off
$\Delta_{1,n}=\bigcup
\limits_{i=0}^{L_n-1}f^{-i}(B_{\delta_{1,n}}(z_1))$ and there
exists $j_{2,n},j_{3,n}$ such that
\begin{itemize}
 \item[E1] $(g_n^*)^j(b_n)=f^j(b_n)$ for $0\leq j \leq j_{1,n}^*$,

\item[E2] $(g_n^*)^{j_{2,n}}(b_n)\in B_{\delta_{1,n}}(z_1)$,
$(g_n^*)^{j_{3,n}}(b_n)\in B_{\delta_{0,n}/\rho_n}(z_0)$,

\item[E3] $(g_n^*)^j(b_n)\in \bigcup \limits_{i=1}^{N}\Phi_i$ for
$0\leq j\leq j_{2,n}$ and $j_{3,n}-j_{2,n}<i_{0,n}$.
\end{itemize}

\begin{rk}\label{7.8}
Above argument shows that
$\#\{\{(g^*_n)^j(b_n)\}^{j_{3,n}}_{j=0}\bigcap\ (\bigcup
\limits_{i=1}^{N}\Phi_i)^c\}<i_{0,n}$.
\end{rk}

Now we'll use connecting lemma in the neighborhood of $z_0$, let's
consider $f^{j_{1,n}^*}(b_n)$, it's near $x_0$, we know that the
positive $g_n^*$-orbit of $f^{j_{1,n}^*}(b_n)$ hits
$B_{\delta_{0,n}/\rho_n}(z_0)$ after $f^{j_{1,n}^*}(b_n)$ and the
negative $g_n^*$-orbit of $f^{j_{1,n}^*}(b_n)$ hits
$B_{\delta_{0,n}/\rho_n}(z_0)$ also, by connecting lemma, the fact
$\Delta_{0,n}=\bigcup f^j(B_{\delta^\prime_n}(z_0))\subset \bigcup
\limits_{i=1}^N \Phi_i$ and remark \ref{7.4}, there exists $g_n\in
B_{\varepsilon_n}(f)$ such that $g_n\equiv f$ off $\Delta_{0,n}$
and there exists $j_0,j_1$ such that
\begin{itemize}
 \item[F1] $g_n^{j_1}(f^{j_{1,n}^*}(b_n))=g_n^{-j_0}(f^{j_{1,n}^*}(b_n))\in B_{\delta^\prime_n}(z_0)$,

 \item[F2] $f^{j_{1,n}^*-j}(b_n)=(g_n^*)^{-j}(f^{j_{1,n}^*}(b_n))=(g_n)^{-j}(f^{j_{1,n}^*}(b_n))$ for $0\leq j\leq
 j_{1,n}^*-j_{0,n}^*$, it means that \\$\#\{Orb_{g_n}(f^{j_{1,n}^*}(b_n))\bigcap \bigcup
\limits_{i=1}^N \Phi_i\}\geq j_{1,n}^*-j_{0,n}^*$,

 \item[F3] $\#\{Orb(f^{j_{1,n}^*}(b_n))\bigcap (\bigcup
\limits_{i=1}^N \Phi_i)^c\}\leq j_{3,n}-j_{2,n}\leq i_{0,n}$.
\end{itemize}

We denote the above periodic orbits for $g_n$ by $Orb(p_n)$ where
$p_n=g_n^{j_1}(f^{j_{1,n}^*}(b_n))$, so we know that $\lim
\limits_{n\rightarrow \infty}p_n\longrightarrow z_0$ and
$\frac{\#\{Orb_{g_n}(p_n)\bigcap \bigcup \limits_{i=1}^N
\Phi_i\}}{Orb_{g_n}(p_n)}= 1-\frac{\#\{Orb_{g_n}(p_n)\bigcap
(\bigcup \limits_{i=1}^N \Phi_i)^c\}}{Orb_{g_n}(p_n)}\geq
1-\frac{i_{0,n}}{j_{1,n}^*-j_{0,n}^*}>1-(1-s)=s$ .

Now we know that there exists a family of diffeomorphisms
$\{g_n\}$ such that $g_n\overset{C^1}{\longrightarrow}f$ and $g_n$
has periodic point $p_n$ such that
$\frac{\#\{Orb_{g_n}(p_n)\bigcap \bigcup \limits_{i=1}^N
\Phi_i\}}{Orb_{g_n}(p_n)}>s$ and $p_n\longrightarrow z_0$, recall
that $z_0=f^{i_0}(y)$, we know that when $n$ is big enough,
$Orb_{g_n}(p_n)$ will pass through $U_0$ the neighborhood of $y$,
so by generic property lemma \ref{7.2}, $f$ itself has periodic
point $p$ such that $\frac{\#\{Orb(p)\bigcap \bigcup
\limits_{i=1}^N \Phi_i\}}{Orb(p)}>s$ and $Orb(p)\bigcap U_0\neq
\phi$. \qed

\bt{99}
\bib{1ABC}
F. Abdenur, C. Bonatti and S. Crovisier, Global dominated
splittings and the $C^1$ Newhouse phenomenon, {\it Proceedings of
the American Mathematical Society} {\bf 134}, (2006), 2229-2237.\
\bib{1ABCD}
F. Abdenur, C. Bonatti and S. Crovisier, L.J. Diaz and L. Wen,
{\it Periodic points and homoclinic classes,} {\it preprint}
(2006).\
\bib{1AS}
R. Abraham and S. Smale, Nongenericity of $\Omega$ -stability,
Global analysis I, {\it Proc. Symp. Pure Math. AMS} {\bf 14}
(1970), 5-8.\
\bib{1A }
M-C. Arnaud, Creation de connexions en topologie $C^1$. {\it
Ergodic Theory and Dynamical System}{\bf 31} (2001), 339-381.\
\bib{1BC}
C. Bonatti and S. Crovisier, Recurrence et genericite(French),
{\it Invent. math.,} {\bf 158} (2004), 33-104\
\bib{1BDP}
C. Bonatti, L. J. D\'{i}az and E. Pujals, A $C^1$-generic
dichotomy for diffeomorphisms: weak forms of hyperbolicity or
infinitely many class or sources, {\it Ann. of Math.} {\bf 158}
(2003), 355-418\
\bib{1BD}
C. Bonatti, L. J. D\'{i}az, Connexions h\'{e}t\'{e}rocliniques et
g\'{e}n\'{e}ricite d'une infinit\'{e} de puits ou de sources, {\it
Annales Scientifiques de l'¨¦cole Normal Sup¨¦rieure de Paris},
{\bf 32 (4)}, (1999) 135-150,\
\bib{1BDV}
C. Bonatti, L. J. D\'{i}az and M. Viana, Dyanamics beyond uniform
hyperbolic, {\bf Volume 102} of {\it Encyclopaedia of Mathematical
Sciences.} Springer- Verlag, Berlin, 2005. A global geometric and
probabilistic perspective, Mathematical Physics, III.\
\bib{1BGW}
C. Bonatti, S. Gan and L. Wen, On the existence of non-trivial
homoclinic class, {\it preprint} (2005)\
\bib{1BGV}
C. Bonatti, N. Gourmelon and T. Vivier, Perturbation of the
derivative along periodic orbits, {\it preprint} (2004)\
\bib{1Co}
C. Conley, Isolated invariant sets and Morse index, CBM Regional
Conference Series in Mathematics, {\bf 38}, AMS Providence,
R.I.,(1978).\
\bib{1C1}
S. Crovisier, Periodic orbits and chain transitive sets of $C^1$
diffeomorphisms, {\it preprint} (2004).\
\bib{1C2}
S. Crovisier, Birth of homoclinic intersections: a model for the
central dynamics of partial hyperbolic systems, {\it preprint}
(2006).\
\bib{1F}
J. Franks, Necessary conditions for stability of diffeomorphisms,
{\it Trans. Amer. Math. Soc.}, {\bf 158} (1977), 301-308.\
\bib{1G1}
S. Gan, Private talk.\
\bib{1G2}
S. Gan, Another proof for $C^1$ stability conjecture for flows,
{\it SCIENCE IN CHINA (Series A)} {\bf 41 No. 10} (October 1998)
1076-1082 \
\bib{1G3}
S. Gan, The Star Systems $X^*$ and a Proof of the $C^1
Omega$-stability Conjecture for Flows, {\it Journal of
Differential Equations}, {\bf 163} (2000) 1--17\
\bib{1GW1}
S. Gan and L. Wen, Heteroclinic cycles and homoclinic closures for
generic diffeomorphisms, {\it Journal of Dynamics and Differential
Equations,}, {\bf 15} (2003), 451-471.\
\bib{1GW2}
S. Gan and L. Wen, Nonsingular star flow satisfy Axion M and the
nocycle condition, {\it Ivent. Math.,} {\bf 164} (2006), 279-315.\
\bib{1H}
S. Hayashi, Connecting invariant manifolds and the solution of the
$C^1$ stability and $\Omega$-stable conjecture for flows, {\it
Ann. math.}, {\bf 145} (1997), 81-137.\
\bib{1HPS}
M. Hirsch, C. Pugh and M. Shub, Invariant manifolds, volume 583 of
{\it Lect. Notes in Math}. Springer Verlag, New york, 1977\
\bib{1L1}
Shantao Liao, Obstruction sets I, {\it Acta Math. Sinica}, {\bf
23} (1980), 411- 453.\
\bib{1L2}
Shantao Liao, Obstruction sets II, {\it Acta Sci. Natur. Univ.
Pekinensis}, {\bf 2} (1981), 1-36.\
\bib{1L3}
Shantao Liao, On the stability conjecture, {\it Chinese Annals of
Math.}, {\bf 1} (1980), 9-30.(in English)\
\bib{1L4}
Shantao Liao, An existence theorem for periodic orbits, {\it Acta
Sci. Natur. Univ. Pekinensis}, {\bf 1} (1979), 1-20.\
\bib{1L5}
Shantao Liao, Qualitative Theory of Differentiable Dynamical
Systems, {\it China Science Press}, (1996).(in English)\
\bib{1M1}
R. Ma\~{n}\'{e}, Quasi-Anosov diffeomorphisms and hyperbolic
manifolds, {\it Trans. Amer. Math. Soc.}, {\bf 229} (1977),
351-370.\
\bib{1M2}
R. Ma\~{n}\'{e}, Contributions to the stability conjecture, {\it
Topology}, {\bf 17} (1978), 383-396.\
\bib{1M3}
R. Ma\~{n}\'{e}, An ergodic closing lemma, {\it Ann. Math.}, {\bf
116} (1982), 503-540.\
\bib{1M4}
R. Ma\~{n}\'{e}, A proof of the $C^1$ stability conjecture. {\it
Inst. Hautes Etudes Sci. Publ. Math.} {\bf 66} (1988), 161-210.\
\bib{1MP}
C. Morales and M. J. Pacifico, Lyapunov stability of
$\omega$-limit sets, {\it Discrete Contin. Dyn. Syst.} {\bf 8}
(2002), no. 3, 671-674.\
\bib{1N1}
S. Newhouse, Non-density of Axiom A(a) on S 2. {\it Proc. A. M. S.
Symp pure math}, {\bf 14} (1970), 191-202, 335-347.\
\bib{1N2}
S. Newhouse, Diffeomorphisms with infinitely many sinks, {\it
Topology}, {\bf 13}, 9-18, (1974).\
\bib{1PV}
J. Palis and M. Viana, High dimension diffeomorphisms displaying
infinitely sinks, {\it Ann. Math.}, {\bf 140} (1994), 1-71.\
\bib{1Pl}
V. Pliss, On a conjecture due to Smale, {\it Diff. Uravnenija.},
{\bf 8} (1972), 268-282.\
\bib{1P}
C. Pugh, The closing lemma, {\it Amer. J. Math.}, {\bf 89} (1967),
956-1009.\
\bib{1PS1}
E. Pujals and M. Sambarino, Homoclinic tangencies and
hyperbolicity for surface diffeomorphisms, {\it Ann. Math.}, {\bf
151} (2000), 961-1023.\
\bib{1PS2}
E. Pujals and M. Sambarino, Density of hyperbolicity and
tangencies in sectional dissipative regions, {\it preprint}
(2005).\
\bib{1S}
J. Selgrade, Isolated invariant sets for flows on vector bundles,
{\it Trans. Amer. Math. Soc.}, {\bf 203} (1975), 259-390.\
\bib{1SS}
R. Sacker and G. Sell, Existence of dichotomies and invariant
splittings for linear differential systems, {\it J. Diff. Eq.},
{\bf 22} (1976), 478-496.\
\bib{1Sh}
M. Shub, Topological transitive diffeomorphisms in $T^4$, Lecture
Notes in Math. Vol. {\bf 206}, Springer Verlag, 1971.\
\bib{1w0}
L. Wen, On the $C^1$-stability conjecture for flows. {\it Journal
of Differential Equations,} {\bf 129}(1995) 334-357.\
\bib{1W1}
L. Wen, Homoclinic tangencies and dominated splittings, {\it
Nonlinearity}, {\bf 15} (2002), 1445-1469.\
\bib{1W2}
L. Wen, A uniform $C^1$ connecting lemma, {\it Discrete and
continuous dynamical systems}, {\bf 8} (2002), 257-265.\
\bib{1W3}
L. Wen, Generic diffeomorphisms away from homoclinic tangencies
and heterodimensional cycles, {\it Bull. Braz. Math. Soc. (N.S.)},
{\bf 35} (2004), 419- 452.\
\bib{1W4}
L. Wen, Selection of Quasi-hyperbolic strings, {\it preprint}
(2006).\
\bib{1WG}
L. Wen and S. Gan, Obstruction sets, Obstruction sets,
quasihyperbolicity and linear transversality.(Chinese) {\it
Beijing Daxue Xuebao Ziran Kexue Ban}, {\bf 42} (2006), 1-10.\
\bib{1XW}
Z. Xia and L. Wen, $C^1$ connecting lemmas, {\it Trans. Amer.
Math. Soc.} {\bf 352} (2000), 5213-5230.\
\bib{1YGW}
D. Yang, S. Gan, L, Wen, Minimal Non-hyperbolicity and
Index-Completeness, {\it preprine} (2007).\
\bib{1Y1}
J. Yang, Ergodic measure far away from tangency, {\it preprint}
(2007).\
\bib{1Y2}
J. Yang, Newhouse phenomena and homoclinic class, {\it preprint}
(2007).\
\bib{1Y4}
J. Yang, Aperiodic class, {\it preprint} (2007).\
\bib{1ZG}
Y. Zhang and S. Gan, On Ma\~{n}\'{e}'s Proof of the $C^1$
Stability Conjecture, {\it Acta mathematica Sinica, English
Series} Vol. {\bf 21}, No. {\bf 3}, June, 2005, 533-540.\
 \et

\end{document}